\documentclass[journal,twoside,web]{ieeecolor}
\usepackage{generic}
\usepackage{cite}
\usepackage{amsmath,amssymb,amsfonts}
\usepackage{bbm}
\usepackage{mathrsfs}
\usepackage{algorithmic}
\usepackage{graphicx}
\usepackage{multicol}
\usepackage{multirow}
\usepackage{epstopdf}
\usepackage{makecell}
\usepackage{textcomp}
\usepackage{booktabs}

\newtheorem{lemma}{Lemma}

\newtheorem{thm}{Theorem}
\newtheorem{defi}{Definition}

\newtheorem{rem}{Remark}

\def\BibTeX{{\rm B\kern-.05em{\sc i\kern-.025em b}\kern-.08em
    T\kern-.1667em\lower.7ex\hbox{E}\kern-.125emX}}
\begin{document}
\title{From Generalized Gauss Bounds to Distributionally Robust Fault Detection with Unimodality Information}
\author{Chao Shang, \IEEEmembership{Member, IEEE}, Hao Ye, Dexian Huang, and Steven X. Ding
\thanks{This work was supported in part by National Science and Technology Innovation 2030 Major Project of the Ministry of Science and Technology of China under Grant 2018AAA0101604, and National Natural Science Foundation of China under Grant 62003187 and Grant 61873142.}
\thanks{C. Shang, H. Ye, and D. Huang are with Department of Automation, Beijing National Research Center for Information Science and Technology, Tsinghua University, Beijing 100084, China (e-mail: c-shang@tsinghua.edu.cn, haoye@tsinghua.edu.cn, huangdx@tsinghua.edu.cn). }
\thanks{S. X. Ding is with Institute for Automatic Control and Complex Systems (AKS), University of Duisburg-Essen, Bismarckstrasse 81 BB, 47057 Duisburg, Germany (e-mail: steven.ding@uni-due.de)}}

\maketitle

\begin{abstract}
Probabilistic methods have attracted much interest in fault detection design, but its need for complete distributional knowledge is seldomly fulfilled. This has spurred endeavors in distributionally robust fault detection (DRFD) design, which secures robustness against inexact distributions by using moment-based ambiguity sets as a prime modelling tool. However, with the worst-case distribution being implausibly discrete, the resulting design suffers from over-pessimisim and can mask the true fault. This paper aims at developing a new DRFD design scheme with reduced conservatism, by assuming unimodality of the true distribution, a property commonly encountered in real-life practice. To tackle the chance constraint on false alarms, we first attain a new generalized Gauss bound on the probability outside an ellipsoid, which is less conservative than known Chebyshev bounds. As a result, analytical solutions to DRFD design problems are obtained, which are less conservative than known ones disregarding unimodality. We further encode bounded support information into ambiguity sets, derive a tightened multivariate Gauss bound, and develop approximate reformulations of design problems as convex programs. Moreover, the derived generalized Gauss bounds are broadly applicable to versatile change detection tasks for setting alarm thresholds. Results on a laborotary system shown that, the incorporation of unimodality information helps reducing conservatism of distributionally robust design and leads to a better tradeoff between robustness and sensitivity.
\end{abstract}

\begin{IEEEkeywords}
Fault detection, uncertain systems, optimization, unimodality
\end{IEEEkeywords}

\section{Introduction}
\label{sec:introduction}
\IEEEPARstart{W}{ith} the rapidly growing complexity of modern technical systems, requirements of operational safety and reliability in an uncertain environment are becoming more critical. This has stimulated the tremendous development and successful applications of fault detection, diagnosis and fault-tolerant control techniques over the past few decades \cite{dingadvanced2021,blanke2006diagnosis}. From a unified viewpoint, residual generation and residual evaluation are centerpieces of fault detection design \cite{dingadvanced2021}. The former aims to construct an indicator that can sensitively unveil the occurrence of anomalous events in dynamical systems, while the latter decides whether an alarm shall be raised via some detection logic.

The ubiquity of uncertainties raises significant challenges for fault detection design. In an uncertain environment, the false alarm rate (FAR) and fault detection rate (FDR) are indices of immediate interest for performance evaluation under fault-free and faulty conditions. Since FAR and FDR are substaintially probabilities, it is rational to formulate fault detection design problems probabilistically \cite{esfahani2015tractable}, where a detailed description of underlying distribution shall be available. Current endeavors in this vein widely hinge on statistical inference under the Gaussian assumption, e.g. the generalized likelihood ratio test (GLRT), whereby the $\chi^2$-distribution has been used for thresholding; however, in practice the Gaussian assumption itself may be unjustifiable and thus vulnerable. Once the true distribution deviates from the assumed normality, the detection performance can significantly degrade. Specifically, a high FAR can raise the ``alarm flood" issue, eventually leading to mistrust of the alarm system and fatal vulnerability to abnormal events \cite{wang2015overview}. Along an alternative route is the set-membership technique, which accounts for all admissible uncertainty realizations within a norm-bounded set (e.g. zonotope) and then optimizes the worst-case performance \cite{fagarasan2004causal,ingimundarson2009robust}. Despite its distribution-free nature, the ensuing fault detection design may be over-pessimistic due to the absence of necessary statistical information.

The above limitations are being recognized and addressed by distributionally robust optimization (DRO), an emerging roadmap in operations research community \cite{delage2010distributionally,wiesemann2014distributionally,cherukuri2019cooperative,li2020data}. As an intermediate to aforesaid two mainstreams, DRO shows wider applicability in practical situations where only partial stochastic information is available. The crux is to construct an \textit{ambiguity set} as a collection of admissible distributions sharing some common properties such as the moments and support, based on which the worst-case performance is optimized. In this way, the resultant decision can hedge against the ambiguity in probability distributions of unknowns. Moreover, for a large class of DRO problems the convexity of problems can be recaptured, which secures computational tractability \cite{wiesemann2014distributionally}. Such a new uncertainty characterization has also been popularized in systems and control, see e.g. \cite{van2015distributionally,yang2018dynamic,boskos2020data}. Recently, distributionally robust fault detection (DRFD) has been extensively investigated \cite{renganathan2020distributionally,shang2021distributionally,xue2020optimal}, where robust integrated design of residual generator and alarm threshold are obtained, reliably enforcing constraints on FAR, FDR and other indices irrespective of imprecisely known distributions \cite{shang2021distributionally,xue2020optimal,song2020parity}.

The robustness level of DRFD design relies heavily on the ambiguity set, which is mostly constructed based on the mean and covariance in prior work. Such a description caters to the broad interest in using the first two moments to characterize a distribution; however, the induced design always shows over-pessimism, as can be evidenced from \cite{shang2021distributionally} where the true FAR tends to be excessively lower than the tolerance but compromises the sensitivity against faults. In fact, the moment-based ambiguity set encompasses an excessively large class of probability distributions, among which the worst-case one is found to be pathologically discrete; see e.g. \cite{delage2010distributionally,zymler2013distributionally}. In real-world systems, however, it is unlikely that uncertainties, especially disturbances, are governed by discrete distributions with few atoms. Thereby, the robustness against such implausibility is deemed as a main cause for over-conservatism.

Thus, this paper is oriented towards a new distributionally robust design scheme, which alleviates the conservatism and strikes a sensible tradeoff between FAR and FDR. The idea is to integrate the \textit{unimodality} of distributions with moment and support information in the ambiguity set, effectively ruling out unrealistic discrete distributions. Notably, the unimodality assumption implies that \textit{smaller deviations are always likely than larger ones}, an intrinsic property of many known distributions in probability theory. More importantly, unimodality can be evidenced in largely many practical scenarios, e.g. from histograms or scatter plots, which also inherently justifies the extensive usage of Gaussian distributions as an approximation. As such, departing from usual methods based on either Gaussian or norm-bounded assumptions, the unimodality-induced ambiguity set yields a ``coarse" yet practically sound description to uncertainty governed by unknown distributions. An integrated design problem is then formulated with various types of prior knowledge exploited, which maximizes the overall fault detectability while robustly regulating false alarms under inexact and even varied distributions.

We remark that unimodality has been adopted for uncertainty description in diverse fields including DRO \cite{van2016generalized,Li2017Ambiguous}, control theory \cite{barmish1997uniform,lagoa2003probabilistic}, power system operations \cite{summers2015stochastic,Li2019Distributionally,pourahmadi2021distributionally}, and statistics \cite{bertsimas2005optimal,vandenberghe2007generalized,van2019distributionally}. Despite these efforts, resolving the DRFD design problem remains a challenge, which arises largely from evaluating the worst-case FAR over the ambiguity set. In fact, evaluating the FAR is simply \textit{a quantification problem of tail probability} that a random vector deviates from its mean in the multi-dimensional setting. Thus, given the mean and covariance, the worst-case FAR is nothing but an extension of the classic univariate Chebyshev bound \cite{renganathan2020distributionally}. With unimodality further considered, the worst-case FAR can be also viewed as generalizing the Gauss bound. Nevertheless, multivariate unimodal Gauss bounds developed in \cite{van2016generalized,van2019distributionally} are not applicable to fault detection design. The crux is that quadratic evaluation functions are typically used in fault detection and thus the tail probability outside an ellipsoidal region is of interest. However, to the best of our knowledge, until now no results have regarded the multivariate Gauss bound of an ellipsoidal acceptance region, while only the polyhedral case was addressed in previous contributions \cite{van2016generalized,van2019distributionally}.

To lay the groundwork for distributionally robust design, a new multivariate Gauss bound is developed in closed form, which turns out to be strictly tighter than the classic Chebyshev bound without assuming unimodality. Based on this, we attain analytical expressions of solutions to DRFD problems, which strictly improve upon known DRFD schemes disregarding unimodality \cite{shang2021distributionally}. To further reduce conservatism, we then embark on the design problem where uncertainty is additionally known to be bounded. As such, moments, unimodality and support information can be integrally encoded in the strengthened ambiguity set. Based on this, a new tightened multivariate Gauss bound is developed, by solving a tractable convex program. Despite its suboptimality, there always exists a suitable tuning parameter rendering the Gauss bound no higher than that assuming unbounded support. Thereby, the utilization of more distributional information renders the DRFD design more sensitive to anomaly while safely keeping FAR below an acceptable level. Furthermore, it is shown that the applicability of the generalized Gauss bounds goes far beyond fault detection design, in that they are also useful for reliably setting alarm thresholds in versatile change detection tasks in systems and control, such as attack detection in cyber-physical systems (CPS) as well as control performance monitoring. The effectiveness of the developed DRFD schemes is illustrated on a realistic three-tank apparatus.

This paper proceeds as follows. Section II revisits preliminary knowledge of residual generation, unimodality, and Chebyshev/Gauss bounds in probability theory. In Section III, the main results of this paper are presented, while in Section IV case studies are reported. Section V concludes the paper.

\textit{Notation:} Given an integer $s>0$, the augmented vector is defined as $\xi_s(k) = [\xi(k)^\top ~ \xi(k-1)^\top ~\cdots~ \xi(k-s)^\top]^\top$. For a matrix $A$, its null space and Moore-Penrose inverse are denoted by ${\rm ker}(A)$ and $A^\dagger$, respectively, and for a symmetric $A$ its positive semi-definiteness is indicated by $A \succeq 0$. $I_n$ denotes the identity matrix of size $n \times n$. The line segment connecting two points $x,y\in\mathbb{R}^n$ is denoted by $[x,y] \subseteq \mathbb{R}^n$.

\section{Preliminaries}
\subsection{DRFD design perspective}
Let us consider the following linear discrete-time stochastic system:
\begin{equation}
\left \{
\begin{split}
& x(k+1) = Ax(k) + Bu(k) + B_d d(k) + B_f f(k)  \\
& y(k) = Cx(k) + Du(k) + D_d d(k) + D_f f(k)
\end{split}
\right .
\label{eq:LTI}
\end{equation}
where $x \in \mathbb{R}^{n_x}$, $y \in \mathbb{R}^{n_y}$, $u \in \mathbb{R}^{n_u}$, $d \in \mathbb{R}^{n_d}$, and $f \in \mathbb{R}^{n_f}$ stand for the process state, measured output, control input, unknown stochastic disturbance and faults, respectively. System matrices in \eqref{eq:LTI} are assumed to be known and have appropriate dimensions. Meanwhile, $(C,A)$ is observable. To construct a residual generator, a unified way is to adopt the stable kernel representation (SKR) $\mathcal{K}(z)$ based on analytical redundancy \cite{ding2014data}, such that in the absence of faults and disturbances, viz. $f(k) = 0$ and $w(k) = 0$, one obtains
\begin{equation}
\mathcal{K}(z) \begin{bmatrix}
u(z) \\
y(z) \\
\end{bmatrix} \equiv 0,
\end{equation}
where $z$ denotes the time-shift operator. Thus, the residual signal $r$ can be generated as an information carrier that is sensitive to anomalies:
\begin{equation}
r(z) = \mathcal{K}(z) \begin{bmatrix}
u(z) \\
y(z) \\
\end{bmatrix} \in \mathbb{R}^{n_r}.
\label{eq:resGen}
\end{equation}
Given $r$, an alarm will be declared signifying an ongoing anomalous situation once the value of $J(r)$ goes beyond a decision threshold $J_{\rm th}$, i.e. $J(r) > J_{\rm th}$. The quadratic function $J(r) = \| r \|^2$ have been mostly used as a summary statistic for residual evaluation. Parallel to the residual generator \eqref{eq:resGen} is its \textit{design form}, which delineates the law governing the dynamics of residuals:
\begin{equation}
r(k) = P[W d_s(k) + V f_s(k)],
\label{eq:designform}
\end{equation}
where $s>n_x$ is the given order of augmented vectors. For convenience, we denote by $\xi := d_s(k) \in \mathbb{R}^{n}$ the uncertainty following an unknown distribution $\mathbb{P}$, where $n = n_d(s+1)$. To determine coefficient matrices $W$ and $V$, a variety of options are available, e.g. the parity space method \cite{chow1984analytical}, observer-based method \cite{jones1973failure}, and subspace identification \cite{ding2014data}. For a more compelete summary readers are referred to \cite{ding2014data}. They crux of residual generation lies in the derivation of the design matrix $P$ so that residual $r$ manifests both sensitivity to fault $f$ and robustness against disturbance $\xi$, the latter of which can be quantitatively assessed by FAR under routine ``healthy" conditions.

\begin{defi}[FAR]
Given the threshold $J_{\rm th}$, the FAR of the residual generator is defined as
\begin{equation}
{\rm FAR} = \mathbb{P}_{\xi} \left \{ \| r \|^2  > J_{\rm th} | f = 0 \right \}.
\end{equation}
\label{defi:1}
\end{defi}

Due to $r = PW\xi$ in fault-free cases, nuisance false alarms are raised once $\xi$ falls outside an ellipsoidal confidence region since $\xi$ enters quadratically into the constraint. The risk of such undesired events is quantified by FAR, whose computation involves high-dimensional integral and thus entails full knowledge about the true distribution $\mathbb{P}_{\xi}$. For the sake of tractability, $\mathbb{P}_{\xi}$ has been generically modeled as Gaussian, which can differ vastly from the ground-truth. As a consequence, there will be serious miscalculation of probabilities, leading to an inaccurate value of $J_{\rm th}$ based on the $\chi^2$-distribution. This inspires the usage of the so-called ambiguity sets for uncertainty description. It includes a class of plausible distributions resembling the true distribution, in the precise sense that they share certain common properties such as the first two moments.
\begin{defi}[Moment-based ambiguity set, \cite{delage2010distributionally}] Given the support $\Xi \subseteq \mathbb{R}^n$, the estimated mean $\mu_0$ and covariance $S_0$, the moment-based ambiguity set is defined as:
\begin{equation*}
\begin{split}
& \mathcal{D}(\gamma_1, \gamma_2, \Xi) \\
= & \left \{ \mathbb{P}(\mathrm{d}\xi) \left | \begin{split}
& \mathbb{P}\{ \xi \in \Xi \} = 1 \\
& (\mathbb{E}_\mathbb{P}\{ \xi \} - \mu_0)^\top S_0^{-1} (\mathbb{E}_\mathbb{P}\{ \xi \} - \mu_0) \le \gamma_1 \\
& \mathbb{E}_\mathbb{P}\{ (\xi - \mu_0)(\xi - \mu_0)^\top \} \preceq \gamma_2 S_0
\end{split}
\right . \right \},
\end{split}
\end{equation*}
where $\gamma_1 \ge 0$ and $\gamma_2 \ge \max \{ \gamma_1, 1\}$ are size parameters. Deviations from the ``nominal" mean are described by an ellipsoid centered at $\mu_0$ whose size can be adjusted by $\gamma_1 > 0$, while $\gamma_2 S_0$ upper-bounds the second-order moment in a semi-definite sense. These parameters can be effectively tuned using the bootstrap strategy \cite{bertsimas2018data}.
\label{defi:2}
\end{defi}

The construction of $\mathcal{D}(\gamma_1, \gamma_2, \Xi)$ hinges solely on mean-covariance information, thereby hedging against ambiguity in high-order statistics. Assume throughout that $\mu_0 = 0$. This is because with $\|r\|^2$ used for detection purpose, it is necessary to perform centering to obtain zero-mean residuals, which better distinguishes nominal variations from anomalies \cite{ding2014data}. With sufficient samples $\{ \hat{\xi}^{(i)} \}$, one can obtain $S_0 = \sum_{i=1}^N (\hat{\xi}^{(i)} - \mu_0) (\hat{\xi}^{(i)} - \mu_0)^\top / (N - 1)$ in an empirical way. Given $\mathcal{D}$, the DRFD design problem can be formulated as a distributionally robust chance constrained program \cite{shang2021distributionally}:
\begin{equation}
\begin{split}
\max_{P} &\ \rho(P) \\
\mathrm{s.t.} &\ \sup_{\mathbb{P}_{\xi} \in \mathcal{D}}  \mathbb{P}_{\xi} \left \{ \|PW\xi \|^2 > 1 \right \} \le \varepsilon
\end{split}
\label{eq:DRprobdesign0}
\tag{DRFD}
\end{equation}
where $\rho(\cdot)$ is an overall fault detectability metric to be maximized, and $\varepsilon$ is a prescribed upper-bound of FAR, e.g. 0.05, which stands for the highest frequency of false alarms that can be tolerated in engineering practice. This can be decided, for example, by the maximal labor cost that is affordable for effective alarm removal. The threshold $J_{\rm th} = 1$ is trivially adopted because otherwise one could always attain infinitely many tuples $\{P,J_{\rm th}\}$ with identical detection performance. The distributionally robust chance constraint in \eqref{eq:DRprobdesign0} secures that the worst-case FAR over all distributions in $\mathcal{D}$ does not exceed $\varepsilon$. It thus offers a clear control mechanism of FAR for unknown distributions, such that alarm overloading and safety hazard can be reliably circumvented under fault-free conditions. In this way, one strikes a tradeoff between robustness against disturbance and sensitivity against faults. Viable choices of $\rho(\cdot)$ include the Frobenius norm metric \cite{ding2019application}:
\begin{equation}
\rho_1(P) = \| PV \|_F^2 = \mathrm{Tr}\{ V^\top P^\top PV \},
\end{equation}
and the pseudo-determinant metric \cite{shang2021distributionally}:
\begin{equation}
\begin{split}
\rho_2(P) & = \log \mathrm{pdet}(V^\top P^\top PV) \\
& = \log \det(\Lambda^\top U_1^\top P^\top PU_1 \Lambda).
\end{split}
\end{equation}
where $V = U_1 \Lambda U_2^\top$ is the compact singular value decomposition with $\Lambda \in \mathbb{R}^{m_f \times m_f}$ being diagonal and invertible. In a nutshell, using $\rho_1(P)$ and $\rho_2(P)$ amount to maximizing, respectively, the sum of eigenvalues, and the product of positive eigenvalues of $V^\top P^\top PV$. More generally, a weighted combination of both metrics can also be employed.
\begin{rem}
Note that the inexactness of distribution of additive disturbance $\xi = d_s(k)$ is addressed by problem \eqref{eq:DRprobdesign0}. However, it enables a broader usage in general cases where knowledge about $d_s(k)$ itself is unavailable (e.g. unmeasurable disturbance) but a running residual signal $v(k)$ has already been developed, whose dynamics is governed by $v(k) = W_v d_s(k) + V_v f_s(k)$. In this case, the residual in fault-free cases can be viewed as additive uncertainty, viz. $\xi := v_0(k) = W_v d_s(k)$, whose distribution is unknown but data samples can be attained under routine fault-free conditions to construct $\mathcal{D}$. It then follows that $v(k) = \xi(k) + V_v f_s(k)$ and the design goal is to identify a design matrix $P$ to ``refine'' the present residual generator as 
\begin{equation}
r(k) = Pv(k) = P[v_0(k) + V_v f_s(k)].
\label{eq:designform2}
\end{equation}
In this case, problem \eqref{eq:DRprobdesign0} still applies with $W = I$ and $V$ replaced by $V_v$ in the definition of the metric $\rho(P)$.
\label{rem:1}
\end{rem}
\begin{rem}
A large body of work postulate \textit{exact moment matching}, which uses $\mathcal{D}(0, 1, \mathbb{R}^n)$ and replaces ``$\preceq$" with ``$=$" therein. In fact, in the context of DRFD, it suffices to consider $\mathcal{D}(\gamma_1, \gamma_2, \Xi)$ with a semi-definite constraint because in \eqref{eq:DRprobdesign0} the worst-case distribution tends to be maximally spread out, attaining the upper-bound on the covariance. Thus, the results developed under moment ambiguity apply straightforwardly to the setup of exact moment matching. Similar arguments have been made in \cite{delage2010distributionally,zhang2018ambiguous}.
\end{rem}

\subsection{Unimodality of unknown's distributions}
The moment-based ambiguity set is known to be prone to over-conservatism. It was shown in \cite{delage2010distributionally} that the worst-case distribution tends to be ``impulsive", thereby being far away from the reality. In the context of fault detection, this renders FAR much lower than the tolerance while unnecessarily sacrificing fault detectability \cite{shang2021distributionally}. To address this issue, incorporating structural properties such as unimodality and monotonicity has been suggested, see e.g. \cite{van2016generalized,van2019distributionally}. As a minimal structural property, unimodality is not only ubiquitous in real-life situations but also inherited by numerous distributions in probability theory.

Conceptually, a distribution is unimodal if larger deviations are less likely than smaller ones. Next we formalize the definition of unimodality. In the univariate case $\xi \in \mathbb{R}$, unimodality asserts the existence of a mode $m$ where the density culminates, along with the cumulative distribution function (cdf) $F(\xi)$ nondecreasing on $(-\infty, m]$ and nonincreasing on $(m,+\infty)$. In the multivariate case, a straightforward generalization is the star-unimodality, which enforces the density function to be non-increasing along any ray emanating from the mode, thereby being ``bell-shaped" intuitively. A precise definition is made based on the notion of star-shaped sets, thereby allowing for distributions that do not admit density functions.

\begin{defi}[Star-shaped set, \cite{Dharmadhikari1988Unimodality}]
A set $\mathcal{A} \subseteq \mathbb{R}^n$ is called star-shaped with center zero if for every $\xi \in \mathcal{A}$, the line segment $[0,\xi]$ is included in $\mathcal{A}$.
\end{defi}

\begin{defi}[Star-unimodality, \cite{Dharmadhikari1988Unimodality}]
A probability distribution $\mathbb{P}$ is said to be star-unimodal about mode 0 if it belongs to the weak closure of the convex hull of all uniform distributions on zero-centered star-shaped sets.
\end{defi}

For continuous probability distributions, a sufficient and necessary condition of star-unimodality is the non-increasing characteristic of density function along any ray emitted from the mode \cite{Dharmadhikari1988Unimodality}. Thus, the concept of star-unimodality generalizes that of univariate unimodality, based on which an extension of multivariate unimodality can be further developed.

\begin{defi}[$\alpha$-unimodality, \cite{olshen1970generalized}]
For any $\alpha>0$, a multivariate distribution $\mathbb{P}$ is $\alpha$-\textit{unimodal} about $0$ if $t^{\alpha}\mathbb{P}(\mathcal{S}/t)$ is non-decreasing in $t\in(0,\infty)$ for every Borel set $\mathcal{S} \in \mathcal{B}(\mathbb{R}^n)$.
\end{defi}

Beyond the generic star-unimodality, the $\alpha$-unimodality regulates the minimal decreasing rate of density along rays emitted from the mode, interpreted as a characterization of the ``degree of unimodality". When $\alpha = n$, the generic star-unimodality is then recovered. Denoting by $\mathscr{P}_\alpha$ the set of all $\alpha$-unimodal distributions with zero mode, it turns out that $\mathscr{P}_\alpha$ enjoys the nesting property $\mathscr{P}_{\alpha_1} \subseteq \mathscr{P}_{\alpha_2}$ if $\alpha_1 \le \alpha_2 \le \infty$. As $\alpha$ tends to infinity, the restriction on unimodality gets relaxed and eventually vanishes.

In this work, we consider the following \textit{structured ambiguity set} that encodes moment information together with $\alpha$-unimodality:
\begin{equation}
\mathcal{D}_\alpha(\gamma_1, \gamma_2, \Xi) \triangleq  \mathcal{D}(\gamma_1, \gamma_2, \Xi) \cap \mathscr{P}_\alpha.
\end{equation}
Clearly, $\mathcal{D}_\alpha(\gamma_1, \gamma_2, \Xi)$ is less conservative than its unstructured counterpart $\mathcal{D}(\gamma_1, \gamma_2, \Xi)$ for any finite $\alpha$. For clarity, we hereafter refer to $\mathcal{D}_\alpha$ with its dependence on $\{\gamma_1, \gamma_2, \Xi\}$ dropped when no confusion is caused. Note that the dirac distribution $\delta_{a}$ with $a \neq 0$ does not belong to $\mathscr{P}_\alpha$ for any finite $\alpha$ \cite{van2016generalized}. This sheds light on the capability of $\mathcal{D}_\alpha$ in eliminating unrealistic discrete distributions except for $\delta_{0}$. To specify $\alpha$, a ubiquitous case without requiring too much \textit{a priori} knowledge is that disturbance at a higher energy level is less likely. As a consequence, $\mathbb{P}_{\xi}$ is known to be star-unimodal, and thus one can safely choose $\alpha = n$. Moreover, a smaller $\alpha$ can be helpful for further reducing the conservatism of $\mathcal{D}_\alpha$ thanks to the nesting property.

\subsection{Optimal inequalities in probability theory}
By definition, the FAR is essentially a tail probability of unfortunate events. Consider the simple univariate case $\xi \in \mathbb{R}$ with mean $\mu$ and variance $\sigma^2$, the worst-case FAR in \eqref{eq:DRprobdesign0} can be evaluated using the Chebyshev inequality, a fundamental result from probability theory:
\begin{equation}
\mathbb{P} \left \{ |\xi - \mu| \ge \kappa \sigma \right \} \le \min \left \{ {1 \over \kappa^2}, ~1 \right \}.
\label{eq:Chebyshev}
\end{equation}
Due to its distribution-free nature, the Chebyshev inequality constitutes the foundation of a variety of probabilistic methods such as the minimax probability machines and DRFD design \cite{renganathan2020distributionally,xue2020optimal,xue2020distribution}. Note that the Chebyshev bound is tight due to the existence of an extremal distribution making \eqref{eq:Chebyshev} an equality. Such a distribution is known to be discrete, which is intimately related to the interplay between the Chebyshev bound and DRO problems. That is, the r.h.s. of \eqref{eq:Chebyshev} can be viewed as the optimal value of the following worst-case probability problem with the moment-based ambiguity set \cite{van2016generalized}:
\begin{equation}
\sup_{\mathbb{P} \in \mathcal{D}'} \mathbb{P} \left \{ |\xi - \mu| \ge \kappa \sigma \right \},
\label{eq:Chebyshev_DRO}
\end{equation}
where $\mathcal{D}'$ encloses univariate distributions sharing the same mean $\mu$ and variance $\sigma^2$:
\begin{equation}
\mathcal{D}' = \left \{ \mathbb{P}(\mathrm{d}\xi) \left | \mathbb{E}_\mathbb{P} \{ \xi \}= \mu, ~\mathbb{E}_\mathbb{P}\{ (\xi - \mu)^2 \} = \sigma^2
\right . \right \}.
\end{equation}
In 1823, it was first proved by Gauss \cite{gauss1823theoria} that, considering unimodal distributions with the mode coinciding with the mean, the conservatism of the Chebyshev bound can be alleviated by the following Gauss bound
\begin{equation}
\mathbb{P} \left \{ |\xi - \mu| \ge \kappa \sigma \right \} \le \left \{ \begin{split}
& {4 \over 9\kappa^2},~~~~~ {\rm if}~\kappa > {2 \over \sqrt{3}} \\
& 1 - { \kappa \over \sqrt{3}},~  {\rm otherwise}
\end{split} \right .
\label{eq:Gauss}
\end{equation}
with an exact improvement factor $4/9$ for sufficiently large deviations. Similar to its Chebyshev counterpart, the Gauss bound can be interpreted as the worst-case probability problem \eqref{eq:Chebyshev_DRO} with a strengthened $\alpha$-unimodal ambiguity set $\mathcal{D}'_{\alpha} \triangleq \mathcal{D}_{\alpha} \cap \{ \mathbb{P}(\mathrm{d}\xi) | \mathbb{P}~{\rm is}~{\rm unimodal~about~}\mu \}$ in lieu of $\mathcal{D}'$.

There is a vast literature on multivariate generalizations of the Chebyshev bound \cite{bertsimas2005optimal,vandenberghe2007generalized,navarro2014can,navarro2016very}. Meanwhile, multivariate extensions of the univariate Gauss inequality have been investigated as well \cite{meaux1990calculation,van2016generalized,van2019distributionally}, all of which concentrate on the highest risk of a random vector residing outside a polytope. To exemplify, we recall an instance with a closed-form expression.
\begin{thm}\cite[Lemma 4]{van2016generalized}
Consider a hypercube with the center $\mu = [\mu_1,\cdots,\mu_n]^\top$ and the edge length of $2\kappa \sigma$ $(0<\kappa<\infty)$. For multivariate uncertainty $\xi \in \mathbb{R}^n$, the associated $\alpha$-unimodal Gauss bound is explicitly expressed as:
\begin{equation}
\begin{split}
& \left \{ \begin{split}
\sup_{\mathbb{P}} & ~ \mathbb{P} \left \{ \max_{1\le i \le n} |\xi_i - \mu_i| \ge \kappa \sigma \right \} \\
{\rm s.t.} & ~ \mathbb{E}_\mathbb{P} \{ \xi \}= \mu,~ \mathbb{E}_\mathbb{P} \{ \xi\xi^\top \}= {\sigma^2 \over n} I_n + \mu\mu^\top \\
              & ~ \mathbb{P} {\rm ~is~}\alpha{\rm -unimodal~about~}\mu
\end{split} \right . \\
= & \left \{ \begin{split}
& {c_\alpha \over \kappa^2}, ~~~~~~~~~~~~~~~~~~~~~~~{\rm if}~ \kappa > \sqrt{c_\alpha(\alpha+2) \over \alpha}\\
& 1 - \left ( {\alpha \over \alpha + 2} \right )^{\alpha / 2} \kappa^\alpha , ~~{\rm otherwise}
\end{split} \right .
\end{split}
\end{equation}
where the improvement factor is given by
\begin{equation}
c_\alpha = \left ( {2 \over \alpha + 2} \right )^{2 / \alpha}.
\end{equation}
\label{thm:1}
\end{thm}

Theorem \ref{thm:1} subsumes a variety of known probability bounds as special instances. For example, it simplifies to the univariate $\alpha$-unimodal Gauss bound \cite[Theorem 3.3.1]{stellato2014} when $n=1$, and the classic Gauss bound \eqref{eq:Gauss} by setting $\alpha = 1$. As $\alpha \to \infty$, one eventually arrives at the Chebyshev bound \eqref{eq:Chebyshev_DRO} due to
\begin{equation}
\lim_{\alpha \to \infty} c_\alpha = 1,~ \lim_{\alpha \to \infty}  \left ( {\alpha \over \alpha + 2} \right )^{\alpha / 2} \kappa^\alpha = 0.
\end{equation}

\section{Main Results}
\subsection{New multivariate generalization of Gauss inequality}
Insofar as the probability outside an ellipsoid is concerned in the design problem \eqref{eq:DRprobdesign0}, known multivariate $\alpha$-unimodal Gauss bounds no longer apply. To fill this knowledge gap, a new multivariate extension of $\alpha$-unimodal Gauss bounds is first derived in this section. Before proceeding, we recall a useful fact recalled that the family of $\alpha$-unimodal distributions can be reparameterized explicitly based on radial $\alpha$-unimodal distributions as extremal ones.

\begin{defi}[Radial $\alpha$-unimodal distributions, \cite{van2016generalized}]
For any $\alpha > 0$ and $x \in \mathbb{R}^n$, denote by $\delta_{[0,x]}^\alpha$ the radial distribution supported on the line segment $[0,x] \subset \mathbb{R}^n$ with the property $\delta_{[0,x]}^\alpha([0,\lambda x]) = \lambda^\alpha,~ \forall \lambda \in [0,1]$.
\end{defi}


The $\alpha$-unimodality of radial distributions $\delta_{[0,x]}^\alpha(\cdot)$ is rather easy to verify. Moreover, they are extremal distributions in $\mathscr{P}_{\alpha}$ \cite{van2016generalized}, which are not representable as a strict convex combinations of two distinct distributions in $\mathscr{P}_{\alpha}$. Thus, using the Choquet theory \cite{phelps2001lectures}, the family of $\alpha$-unimodal distributions can be explicitly reparameterized by ``mixing" extremal ones.
\begin{lemma}[Choquet representation of unimodal distributions]
For every distribution $\mathbb{P} \in \mathscr{P}_\alpha$ supported on $\Xi$, there exists a unique distribution $\mathbb{Q}$ supported on $\Xi$ such that
\begin{equation}
\mathbb{P}({\rm d} \xi) = \int_{\Xi} \delta_{[0,w]}^\alpha({\rm d} \xi) \mathbb{Q}({\rm d}w).
\label{eq:15}
\end{equation}
\label{thm:2}
\end{lemma}
\begin{proof}
The proof of \cite[Theorem 3.5]{Dharmadhikari1988Unimodality} applies with minor modifications to the present setup and is thus omitted.
\end{proof}

Lemma \ref{thm:2} asserts that every $\alpha$-unimodal distribution supported on $\Xi$ is expressible as a mixture of radial distributions $\delta_{[0,w]}^\alpha(\cdot),~w \in \Xi$, with $\mathbb{Q}$ being the mixture distribution. This allows to recast the worst-case probability problem explicitly in the following supporting lemma.

\begin{lemma}
Given an ellipsoidal confidence region $\mathcal{E} = \{ \xi \in \mathbb{R}^n | \xi^\top M \xi \le 1 \}$ with $M \succeq 0$, the worst-case probability outside $\mathcal{E}$
\begin{equation*}
\sup_{\mathbb{P} \in \mathcal{D}_\alpha(\gamma_1,\gamma_2,\Xi)} \mathbb{P}\{ \xi \notin \mathcal{E} \}
\end{equation*}
is equal to the optimal value of the following semi-infinite optimization problem:
\begin{equation}
\begin{split}
\min_{Q,q,q_0} &\ \gamma_2 {\rm Tr} \left \{
QS^\alpha_0 \right \} + q_0  \\
\mathrm{s.t.}~ &~ \xi^\top Q\xi + 2\xi^\top q + q_0 \ge L_\alpha(\xi) ,~ \forall \xi \in \Xi \\
                       &~ Q \succeq 0
\end{split}
\label{eq:5}
\end{equation}
where $$L_\alpha(\xi) = \max \{1 - \| M^\frac12 \xi \|^{-\alpha}, 0 \},~ S^\alpha_0 = \frac{\alpha + 2}{\alpha} S_0.$$
\label{thm:3}
\end{lemma}

\begin{proof}
We draw ideas from \cite[Theorem 1]{van2019distributionally}. Thanks to Lemma \ref{thm:2}, it suffices to optimize over the unstructured mixture distribution $\mathbb{Q}$ instead of $\mathbb{P}$. Using the reparameterization \eqref{eq:15} yields:
\begin{equation*}
\begin{split}
& \sup_{\mathbb{P} \in \mathcal{D}_\alpha(\gamma_1,\gamma_2,\Xi)} \mathbb{P}_{\xi}  \left \{ \xi^\top M\xi > 1 \right \} \\
= & \sup_{\mathbb{P} \in \mathcal{D}_\alpha(\gamma_1,\gamma_2,\Xi)} \int_{\Xi} \mathbbm{1}_{\{ \xi^\top M\xi > 1 \}}(\xi) \mathbb{P}({\rm d}\xi) \\
= & \left \{ \begin{split}
\sup_{\mathbb{Q}} & \int_{\Xi} \int_{0}^1 \mathbbm{1}_{\{ \xi^\top M\xi > 1 \}}(\lambda w) \delta_{[0,w]}^\alpha(w{\rm d}\lambda) \mathbb{Q}({\rm d}w) \\
{\rm s.t.} & \int_{\Xi} \delta_{[0,w]}^\alpha({\rm d}\xi) \mathbb{Q}({\rm d}w) \in \mathcal{D}(\gamma_1,\gamma_2,\Xi)
\end{split} \right .
\end{split}
\end{equation*}
where the last equality stems from the fact that $\delta_{[0,w]}^\alpha(\xi)$ is supported on the line segment $[0,w]$ only, and thus it suffices to integrate over $[0,w]$ using a scalar $\lambda$. Then the objective amounts to
\begin{equation*}
\begin{split}
& ~\int_{\Xi} \int_{0}^1 \mathbbm{1}_{\{ \xi^\top M\xi > 1 \}}(\lambda w) \delta_{[0,w]}^\alpha(w{\rm d}\lambda) \mathbb{Q}({\rm d}w) \\
= & ~ \int_{\Xi} \max \{1 - \| M^\frac12 \xi \|^{-\alpha}, 0 \} \mathbb{Q}({\rm d}w) \\
\triangleq & ~ \int_{\Xi} L_{\alpha}(w) \mathbb{Q}({\rm d}w) \\
= & ~ \mathbb{E}_{\mathbb{Q}} \left \{ L_{\alpha}(w) \right \}
\end{split}
\end{equation*}
By similar arguments, two moment constraints can be translated into
\begin{equation*}
\begin{split}
& \int_{\Xi} ww^\top \mathbb{Q}({\rm d}w) = \frac{S_0}{\int_0^1 \lambda^2 \delta_{[0,w]}^\alpha(w{\rm d}\lambda)} = \frac{\alpha + 2}{\alpha} S_0 \triangleq S^\alpha_0, \\
& \int_{\Xi} w \mathbb{Q}({\rm d}w) = \frac{\mu_0}{\int_0^1 \lambda \delta_{[0,w]}^\alpha(w{\rm d}\lambda)} = \frac{\alpha + 1}{\alpha} \mu_0 \triangleq \mu_0^\alpha.
\end{split}
\end{equation*}
In this way, the worst-case FAR problem can be recast as the following worst-case expectation problem with the unstructured distribution $\mathbb{Q}$ being the decision variable, which resides in the generic moment-based set:
\begin{equation*}
\begin{split}
\sup_{\mathbb{Q}} &~  \mathbb{E}_{\mathbb{Q}} \left \{ L_{\alpha}(w) \right \} \\
{\rm s.t.} &~ \mathbb{Q} \in \mathcal{D}(S^\alpha_0, \mu_0^\alpha, \Xi)
\end{split}
\end{equation*}
Its full equivalence with problem \eqref{eq:5} immediately follows from the duality argument in \cite[Lemma 1]{delage2010distributionally}, which completes the proof.
\end{proof}

By resolving problem \eqref{eq:5} with unbounded support $\Xi = \mathbb{R}^n$, we arrive at a new multivariate $\alpha$-unimodal Gauss bound in closed form, which not only complements Theorem \ref{thm:1} but also underpins subsequent DRFD design.
\begin{thm}[Generalized $\alpha$-unimodal Gauss bound] Consider $\mathcal{D} = \mathcal{D}_\alpha(\gamma_1,\gamma_2,\mathbb{R}^n)$ that includes all $\alpha$-unimodal distributions on unbounded support $\Xi = \mathbb{R}^n$ subject to moment constraints. The induced worst-case probability of the event $\xi \notin \mathcal{E}$ is upper-bounded by
\begin{equation}
\begin{split}
& ~\sup_{\mathbb{P} \in \mathcal{D}_\alpha(\gamma_1,\gamma_2,\mathbb{R}^n)} \mathbb{P}\{ \xi \notin \mathcal{E} \} \\
\le & \left \{ \begin{split}
& c_\alpha \gamma_2 {\rm Tr} \{ MS_0 \},~~~~~~~~~~{\rm if}~c_\alpha \gamma_2{\rm Tr} \{ MS_0 \} \le {\alpha \over \alpha + 2}  \\
& 1 - \frac{1}{(\gamma_2{\rm Tr} \{ MS_0^{\alpha} \})^{\alpha / 2}},~{\rm otherwise}
\end{split} \right .
\end{split}
\label{eq:10}
\end{equation}
\label{thm:4}
\end{thm}

The proof will be deferred to Appendix due to its complexity. Next we dwell on the effect of introducing $\alpha$-unimodality. In the limiting case $\alpha \to \infty$, where $\alpha$-unimodality vanishes, the following established result is recalled.
\begin{thm}[Generalized Chebyshev bound, \cite{navarro2014can,navarro2016very,shang2021distributionally}]
Consider all distributions $\mathbb{P} \in \mathcal{D}(\gamma_1,\gamma_2,\mathbb{R}^n)$ that have unbounded support and are subject to moment constraints only. The worst-case probability of the event $\xi \notin \mathcal{E}$ is given by:
\begin{equation}
\sup_{\mathbb{P} \in \mathcal{D}(\gamma_1,\gamma_2,\mathbb{R}^n)} \mathbb{P}\{ \xi \notin \mathcal{E} \} = \min \{ \gamma_2\mathrm{Tr}\{M S_0 \} ,1 \}.
\label{eq:30}
\end{equation}
\label{thm:5}
\end{thm}

\begin{rem} Theorem \ref{thm:5} offers an equality, whereas under $\alpha$-unimodality assumption one merely attains via Theorem \ref{thm:4} an upper-bound that is not necessarily sharp. One may wonder whether or not restricting the search within unimodal distributions leads to a more conservative quantification of tail probability. Interestingly, even though its tightness is elusive at present, the generalized Gauss bound \eqref{eq:10} is \textit{provably lower than} its Chebyshev counterpart \eqref{eq:30} for any finite $\alpha$. To see this, it follows from $c_\alpha < 1$ that $c_\alpha \gamma_2{\rm Tr} \{ MS_0 \} < \gamma_2\mathrm{Tr}\{M S_0 \}$. Then two cases are distinguished. When $c_\alpha \gamma_2{\rm Tr} \{ MS_0 \} \le {\alpha \over \alpha + 2} < 1$, one obtains:
\begin{equation*}
c_\alpha \gamma_2{\rm Tr} \{ MS_0 \} < \min \{ \gamma_2\mathrm{Tr}\{M S_0 \} ,1 \}.
\end{equation*}
As for the case $c_\alpha \gamma_2{\rm Tr} \{ MS_0 \} > {\alpha \over \alpha + 2}$, we have:
\begin{equation*}
\begin{split}
& 1 - \frac{1}{(\gamma_2{\rm Tr} \{ M S_0^\alpha \})^{\alpha / 2}}  \\
< & \min \left \{ c_\alpha \gamma_2 {\rm Tr} \{ M S_0 \} , 1\right \} \\
\le & \min \{ \gamma_2\mathrm{Tr}\{M S_0 \} ,1 \}
\end{split}
\end{equation*}
On the other hand, when $\alpha \to \infty$, one can verify that the upper-bound in \eqref{eq:10} becomes the generalized Chebyshev bound \eqref{eq:30}, which is known to be tight. Thus, we conjecture that the generalized $\alpha$-unimodal Gauss bound \eqref{eq:10} is also tight, which will be deferred to further investigation.
\label{rem:2}
\end{rem}

\begin{rem}
The proposed generalized $\alpha$-unimodal Gauss bound \eqref{eq:10} improves upon its Chebyshev counterpart \eqref{eq:30} by a factor of $c_{\alpha}$ when $\gamma_2{\rm Tr} \{ MS_0 \}$ is suitably small. This coincides with the gain within known multivariate formulations of $\alpha$-unimodal Gauss bounds \cite{van2016generalized,meaux1990calculation}. In Fig. \ref{fig:impr_alpha} we depict the value of $c_\alpha$ under varied $\alpha$, where for a moderately valued $\alpha$ (e.g. between $5$ and $50$), the improvement factor ranges from $0.6$ to $0.8$. Two particular cases are noteworthy. When $\alpha = 1$, the maximal improvement is exactly $4/9$, which closely resembles \eqref{eq:Gauss} and thus can be considered as an extension of the univariate Gauss bound \eqref{eq:Gauss}. When $\alpha \to \infty$, the value of $c_\alpha$ tends to one, and thus the generalized Chebyshev bound \eqref{eq:30} is recovered.
\begin{figure}[ht]
  \centering
  \includegraphics[width=0.4\textwidth]{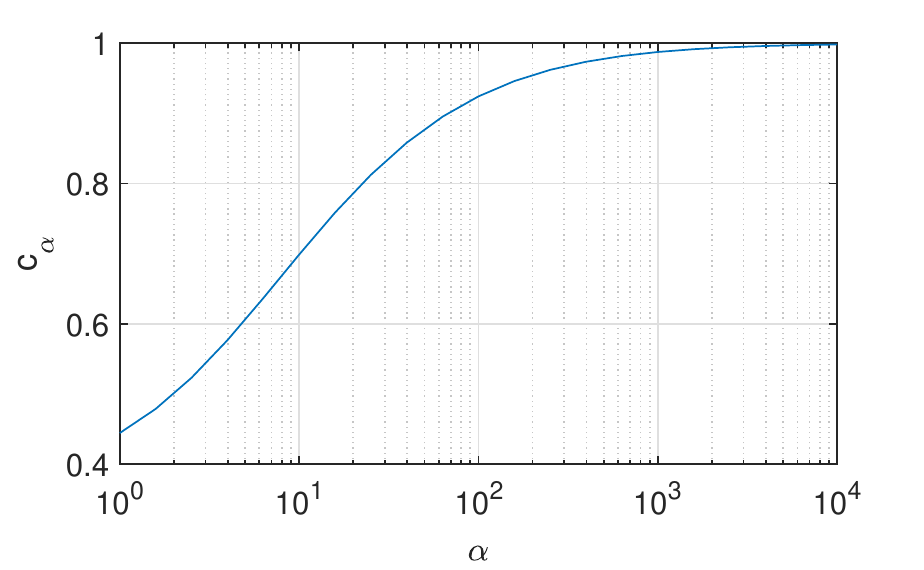}\\
  \caption{Improvement factor $c_\alpha$ of the $\alpha$-unimodal Gauss bound \eqref{eq:10} over its Chebyshev counterpart \eqref{eq:30}.}\label{fig:impr_alpha}
\end{figure}
\end{rem}

\subsection{Solving DRFD problem under unbounded support}
Based on the generalized $\alpha$-unimodal Gauss bound \eqref{eq:10}, we are now in a position to tackle the design problem \eqref{eq:DRprobdesign0} under the structured ambiguity set $\mathcal{D}_\alpha$. It turns out that under different detectability metrics $\rho(\cdot)$, feasible solutions can be derived in closed form, which strictly improve upon known results disregarding unimodality.

\begin{thm}
A feasible solution to problem \eqref{eq:DRprobdesign0} with the Frobenius norm metric $\rho_1(\cdot)$ under $\mathcal{D} = \mathcal{D}_{\alpha}(\gamma_1,\gamma_2,\mathbb{R}^n)$ is given by:
\begin{equation}
P^* = \left \{ \begin{split}
& \sqrt{\omega_1 \varepsilon \over \gamma_2 c_\alpha} \cdot \frac{p_1^\top}{\|V^\top p_1 \|},~~~~~~~~~~~ 0 < \varepsilon \le {\alpha \over \alpha + 2}  \\
& \sqrt{\omega_1 \over \gamma_2 (1 - \varepsilon)^{2/\alpha}} \cdot \frac{p_1^\top}{\|V^\top p_1 \|},~ {\alpha \over \alpha + 2} < \varepsilon \le 1
\end{split} \right .
\label{eq:l2design}
\end{equation}
where $\{ \omega_1,p_1 \}$ are the largest eigenvalue and the associated eigenvector of the generalized eigen-decomposition problem $VV^\top p = \omega W^\top S_0 W^\top p$.
\label{thm:6}
\end{thm}
\begin{proof}
We first consider the case $0 < \varepsilon \le \alpha / (\alpha + 2)$. Noting that the right-hand side of \eqref{eq:10} is strictly increasing in $\gamma_2{\rm Tr}\{ MS_0 \}$ and plugging $M = W^\top P^\top P W$ into \eqref{eq:10}, the constraint on the worst-case FAR in \eqref{eq:DRprobdesign0} is implied by the inequality $c_\alpha \gamma_2 {\rm Tr} \{ W^\top P^\top P W S_0 \} \le \varepsilon.$ Thus, solving the following problem always yields a feasible solution to \eqref{eq:DRprobdesign0} with metric $\rho_1(\cdot)$:
\begin{equation*}
\begin{split}
\max_{P} &\ \mathrm{Tr}\{ V^\top  P^\top P V \} \\
\mathrm{s.t.} &\ c_\alpha \gamma_2 {\rm Tr} \{ W^\top P^\top P W S_0 \} \le \varepsilon
\end{split}
\end{equation*}
One of its optimal solutions is given by the first case in \eqref{eq:l2design} according to \cite[Theorem 2]{shang2021distributionally}. Next we embark on the case $\alpha / (\alpha + 2) < \varepsilon \le 1$. The constraint on the worse-case FAR is implied by $$1 - {1 \over (\gamma_2{\rm Tr} \{ MS_0^{\alpha} \})^{\alpha / 2}} \le \varepsilon,$$ which amounts to $\gamma_2 {\rm Tr} \{ W^\top P^\top P W S_0 \} \le {1 \over (1 - \varepsilon)^{2/\alpha}}.$ By the same token, a feasible solution to the design problem is attained as the second case in \eqref{eq:l2design}, from which the claim follows.
\end{proof}

\begin{thm}
A feasible solution to problem \eqref{eq:DRprobdesign0} with the pseudo-determinant metric $\rho_2(\cdot)$ under $\mathcal{D} = \mathcal{D}_{\alpha}(\gamma_1,\gamma_2,\mathbb{R}^n)$ is given by:
\begin{equation}
P^* = \left \{ \begin{split}  
& \sqrt{\varepsilon \over m_f \gamma_2 c_\alpha} P_{\rm GLRT}, ~~~~~~~~~~~ 0 < \varepsilon \le {\alpha \over \alpha + 2} \\
& \sqrt{1 \over m_f \gamma_2 (1 - \varepsilon)^{2/\alpha}} P_{\rm GLRT},~ {\alpha \over \alpha + 2} < \varepsilon \le 1
\end{split}  \right .
\label{eq:GLRTdesign}
\end{equation}
where $P_{\rm GLRT} = \bar{S}_0^{-1/2}V(V^\top \bar{S}_0^{-1} V)^{\dagger} V^\top \bar{S}_0^{-1}$, $\bar{S}_0 = W S_0 W^\top$ is celebrated GLRT design for anomaly detection \cite{willsky1976generalized,tornqvist2006eliminating}.
\label{thm:7}
\end{thm}
\begin{proof}
We follow the same outline as Theorem \ref{thm:6} and take the case $0 < \varepsilon \le \alpha / (\alpha + 2)$ as an example. In this case, the following problem acts as a conservative approximation to \eqref{eq:DRprobdesign0} with metric $\rho_2(\cdot)$ used:
\begin{equation*}
\begin{split}
\max_{P} &\ \log \det(\Lambda_f^\top U_f^\top P^\top PU_f \Lambda_f) \\
\mathrm{s.t.} &\ c_\alpha \gamma_2 {\rm Tr} \{ W^\top P^\top P W S_0 \} \le \varepsilon
\end{split}
\end{equation*}
which admits an closed-form optimal solution in terms of the first case in \eqref{eq:GLRTdesign} due to \cite[Theorem 3]{shang2021distributionally}. The case $\alpha / (\alpha + 2) < \varepsilon \le 1$ can be treated in a similar manner.
\end{proof}

Some remarks on Theorems \ref{thm:6} and \ref{thm:7} are made in order.
\begin{rem}
When unimodality is disregarded, i.e. $\mathcal{D} = \mathcal{D}(\gamma_1,\gamma_2,\mathbb{R}^n)$, global optimal solutions to problem \eqref{eq:DRprobdesign0} with metrics $\rho_1(\cdot)$ and $\rho_2(\cdot)$ are given by $P^* = \sqrt{\omega_1 \varepsilon / \gamma_2} \cdot p_1^\top / \|V^\top p_1 \|$ and $P^* = \sqrt{\varepsilon / m_f \gamma_2} P_{\rm GLRT}$, respectively \cite[Theorems 2, 3]{shang2021distributionally}. It turns out that they have severer conservatism than the feasible solutions \eqref{eq:l2design} and \eqref{eq:GLRTdesign} derived under $\alpha$-unimodality, in the light of the fact that $0 < c_\alpha < 1$ and $\varepsilon < 1 < 1/(1-\varepsilon)^{2/\alpha}$. When $\varepsilon$ is suitably small, the improvement factor of \eqref{eq:l2design} and \eqref{eq:GLRTdesign} is exactly $c_\alpha$, which can be conceived a ``compensator'' of $\gamma_2$ that describes the ambiguity of covariance estimation and leads to an increased amount of conservatism. 
\end{rem}

\begin{rem}
Two limiting cases are noteworthy. As $\alpha \to 0$, the probability density tends to be concentrated around the mode and collapsed into the dirac distribution $\delta_{ 0 }$. In this case, it is expected that an ``infinitely large'' $P^*$ is attained. This can also be inspected from \eqref{eq:l2design} and \eqref{eq:GLRTdesign} since for $\varepsilon > 0$, ${\alpha \over \alpha + 2} < \varepsilon \le 1$ always holds, resulting in $\lim_{\alpha \to 0}1/(1-\varepsilon)^{2/\alpha} = \infty$. When $\alpha$-unimodality vanishes, i.e. $\alpha \to \infty$, feasible solutions in \eqref{eq:l2design} and \eqref{eq:GLRTdesign} become \textit{global} optimal solutions due to $\lim_{\alpha \to \infty} 1 / \sqrt{c_\alpha} = 1$.
\end{rem}

\begin{rem}
The GLRT design $P_{\rm GLRT}$ itself is developed based on the Gaussian assumption, making $J(r)$ follow a $\chi^2$-distribution with $m_f$ degrees of freedom \cite{willsky1976generalized,tornqvist2006eliminating}. Under confidence level $\varepsilon \in (0,1)$, a theoretical threshold for $J(r)$ is given by $\chi^2_{m_f,1-\varepsilon}$, which enjoys the \textit{strong concentration property} $\chi^2_{m_f,1-\varepsilon} \sim \mathcal{O}(\log(1/\varepsilon))$ as $\varepsilon \to 0$ \cite{laurent2000adaptive}. Under distributional ambiguity, however, such a property no longer remains since Theorem \ref{thm:7} implies that for a sufficiently small $\varepsilon$, a safe threshold for $P_{\rm GLRT}$ is $m_f\gamma_2 c_\alpha / \varepsilon \sim \mathcal{O}(1/\varepsilon)$ as $\varepsilon \to 0$, which is the price we have to pay for being distributionally robust.
\end{rem}


\subsection{DRFD problem under bounded support information}
In practice, stochastic disturbance within a finite time period typically has a bounded energy. This justifies the prevalence of set-membership regime in model-based filtering and fault diagnosis \cite{fagarasan2004causal,ingimundarson2009robust,wang2018zonotopic,ding2019set}, where uncertainty is confined to a compact set. Such knowledge can be further utilized to enrich the information within ambiguity sets by endowing $\mathcal{D}_\alpha$ with a bounded support $\Xi_b \subset \mathbb{R}^n$. By a confluence of moment, unimodality, and bounded support information, the resulting ambiguity set $\mathcal{D}_{\alpha}(\gamma_1,\gamma_2,\Xi_b)$ offers a ``hybrid" description to uncertainty, which can be understood as ``interpolating" between the Gaussian assumption and norm-bounded description. On the one hand, unimodality as well as mean-covariance information underlying Gaussian distributions is preserved. On the other hand, bounded support information is taken into account. All information is seamlessly synthesized by $\mathcal{D}_{\alpha}(\gamma_1,\gamma_2,\Xi_b)$ to robustify the fault detection design while enhancing the detectability. To this end, we first derive a strengthened $\alpha$-unimodal Gauss bound under bounded support, which admits a tractable approximation as a convex program.

\begin{thm}\textit{(Generalized $\alpha$-unimodal Gauss bound with bounded support)}
Suppose the support of distributions is described as the intersection of ellipsoids:
\begin{equation}
\Xi_b = \{ \xi | (\xi - a_j)^\top \Theta_j(\xi - a_j) \le 1,\ j \in \mathbb{N}_{1:N_e} \}.
\label{eq:support}
\end{equation}
Then the worst-case probability of the event $\xi \notin \mathcal{E}$ under $\mathcal{D} = \mathcal{D}_\alpha(\gamma_1,\gamma_2,\Xi_b)$ is no higher than by the optimal value of the following semi-definite program (SDP):
\begin{subequations}
\begin{align}
\min_{Q,q,q_0,\eta,\beta,\tilde{\beta}} & ~ \gamma_2 {\rm Tr} \left \{
QS_0^\alpha \right \} + q_0 \\
\mathrm{s.t.}~~~~ &~ \eta \ge 0,~\beta_j \ge 0,~\tilde{\beta}_j \ge 0,~j=1,\cdots,N_e \label{eq:probdesign16b} \\
\begin{split}
& \begin{bmatrix}
Q - \eta M & q & 0 \\
q^\top & q_0 - 1 & 0 \\
0 & 0 & \eta \\
\end{bmatrix} + \sum_{j=1}^{N_e} \beta_j \Phi_j  \\
& ~~~~~~~~~~~~~~~~ +{1 \over \tau_0^{\alpha}} \begin{bmatrix}
0 & 0 & 0 \\
0 & \alpha + 1 & {-\alpha \over 2\tau_0} \\
0 & {-\alpha \over 2\tau_0} & 0 \\
\end{bmatrix}  \succeq 0 \end{split} \label{eq:probdesign16c} \\
&~ \begin{bmatrix}
Q & q \\
q^\top & q_0
\end{bmatrix} + \sum_{j=1}^{N_e} \tilde{\beta}_j \tilde{\Phi}_j \succeq 0 \label{eq:probdesign16d} \\
&~~ Q \succeq 0 \label{eq:probdesign16e}
\end{align}
\label{eq:probdesign16}
\end{subequations}
where
\begin{equation*}
\tilde{\Phi}_j = \begin{bmatrix} \Theta_j & -\Theta_ja_j \\ -a_j^\top\Theta_j & a_j^\top\Theta_j a_j - 1 \end{bmatrix}, ~\Phi_j = \begin{bmatrix}
\tilde{\Phi}_j & 0 \\
0 & 0 \\
\end{bmatrix}, ~ \tau_0 > 0.
\end{equation*}
\label{thm:8}
\end{thm}
\begin{proof}
By virtue of Lemma \ref{thm:3}, the worst-case probability problem can be recast as:
\begin{equation*}
\begin{split}
\min_{Q,q,q_0} & ~ \gamma_2 {\rm Tr} \left \{
QS_0^\alpha \right \} + q_0 \\
\mathrm{s.t.}~ &\ \xi^\top Q\xi + 2\xi^\top q + q_0 \ge 1 - \| M^\frac12 \xi \|^{-\alpha},~ \forall \xi \in \Xi_b \\
&~ \xi^\top Q\xi + 2\xi^\top q + q_0 \ge 0 ,~ \forall \xi \in \Xi_b \\
&~ Q \succeq 0
\end{split}
\end{equation*}

We first tackle the semi-infinite constraints $\xi^\top Q\xi + 2\xi^\top q + q_0 \ge 1 - \| M^\frac12 \xi \|^{-\alpha},~ \forall \xi \in \Xi_b$, which are implied by the following constraint by the same linearization technique in Theorem \ref{thm:4} (see Appendix):
\begin{equation*}
\xi^\top Q\xi + 2\xi^\top q + q_0 + (\alpha + 1)\tau_0^{-\alpha}  - 1 \ge \alpha \tau_0^{-\alpha - 1} \| M^\frac12 \xi \| ,~ \forall \xi \in \Xi_b
\end{equation*}
for an arbitrary $\tau_0 > 0$. This turns out to be equivalent to the following semi-infinite constraint without conic terms
\begin{equation*}
\begin{split}
& \xi^\top Q\xi + 2\xi^\top q + q_0 + (\alpha + 1)\tau_0^{-\alpha} - \alpha \tau_0^{-\alpha - 1} t - 1 \ge 0 \\
& \forall (\xi, t) \in \left \{ (\xi,t) \left | \begin{split}
& (\xi - a_j)^\top \Theta_j(\xi - a_j) \le 1,\ j \in \mathbb{N}_{1:N_e}  \\
& t^2 \le \xi^\top M\xi \end{split} \right . \right \}. \\
\end{split}
\end{equation*}
By invoking the S-procedure \cite{yakubovich1997s}, a sufficient condition for the above semi-infinite constraint to hold is the existence of Lagrangian multipliers $\{ \beta_j \ge 0,~j=1,\cdots,N_e \}$ and $\eta \ge 0$ such that $\forall (\xi, t) \in \mathbb{R}^{n+1}$,
\begin{equation*}
\begin{split}
& \xi^\top Q\xi + 2\xi^\top q + q_0 - 1 + (\alpha + 1)\tau_0^{-\alpha} - \alpha \tau_0^{-\alpha - 1} t \\
& ~~~~+ \eta t^2 - \eta \xi^\top M\xi + \sum_{j=1}^{N_e} \beta_j (\xi - a_j)^\top \Theta_j (\xi - a_j) \ge 0
\end{split}
\end{equation*}
which amounts to the linear matrix inequality (LMI) \eqref{eq:probdesign16c}.

By similar arguments, the semi-infinite constraints $\xi^\top Q\xi + 2\xi^\top q + q_0 \ge 0,~ \forall \xi \in \Xi_b$ are implied by the existence of $\{ \tilde{\beta}_j \ge 0,~j=1,\cdots,N_e \}$ ensuring \eqref{eq:probdesign16d}. This completes the proof.
\end{proof}

\begin{rem}
It is not restrictive to express $\Xi_b$ as an intersection of ellipsoids. Indeed, along a similar route one can show that the case of a polytopic support $\Xi_b = \{ \xi | G\xi \le h\}$ is also tractable by invoking the S-procedure.
\end{rem}

Recall that in Theorems \ref{thm:4} and \ref{thm:8} upper-bounds of worst-case probabilities are derived, whose tightness is unwarranted. Thus, it is doubtful whether the exploration of support information does help to refine the $\alpha$-unimodal Gauss bound \eqref{eq:10}. Said another way, a possible yet unwanted outcome is that the upper-bound of Theorem \ref{thm:8} is even higher than that of Theorem \ref{thm:4}. Next we show that this case can be effectively circumvented by choosing a suitable $\tau_0 > 0$.
\begin{thm}
Given
\begin{equation}
\tau_0 = \max \left \{ {1 \over \sqrt{c_\alpha}}, \sqrt{\gamma_2{\rm Tr} \{ MS_0^\alpha \}} \right \},
\label{eq:optTau0}
\end{equation}
the optimal value of the associated problem \eqref{eq:probdesign16} is always lower than or equal to the generalized $\alpha$-unimodal Gauss bound \eqref{eq:10} under unbounded support information.
\label{thm:9}
\end{thm}
\begin{proof}
We first consider the case where $\tau_0 = 1 / \sqrt{c_\alpha} > \sqrt{\gamma_2 {\rm Tr} \{ MS_0^\alpha \}}$. From the proof of Theorem \ref{thm:4}, the generalized $\alpha$-unimodal Gauss bound $c_\alpha \gamma_2 {\rm Tr} \{ MS_0 \}$ in \eqref{eq:10} coincides with the optimal value of problem \eqref{eq:thm1_apprx1} with $\tau_0 = 1 / \sqrt{c_\alpha}$. Denote by $\{Q^*,q^*,q_0^*,\eta^* \}$ the associated optimal solution to \eqref{eq:thm1_apprx1}. It immediately follows that $\{Q^*,q^*,q_0^*,\eta^*,\beta_j = 0,\tilde{\beta}_j = 0 \}$ is feasible for \eqref{eq:probdesign16}. Henceforth, \eqref{eq:probdesign16} is less constrained than \eqref{eq:thm1_apprx1} with the same objective, indicating that the optimal value of \eqref{eq:probdesign16} is no lower than that of \eqref{eq:thm1_apprx1}. The case of $\tau_0 = \sqrt{\gamma_2 {\rm Tr} \{ MS_0^\alpha \}} \ge 1 / \sqrt{c_\alpha}$ can be treated in a similar manner, from which the claim follows.
\end{proof}

The usefulness of the SDP formulation \eqref{eq:probdesign16} lies in that it can be modularly embeded in general distributionally robust chance constrained programs to derive tractable approximations. For instance, a tractable approximation of problem \eqref{eq:DRprobdesign0} under $\mathcal{D} = \mathcal{D}_\alpha(\gamma_1,\gamma_2,\Xi_b)$ can be developed, by defining positive semi-definite matrix $\bar{P} = P^\top P \succeq 0$ as the optimization variable:
\begin{equation}
\begin{split}
\max_{\bar{P} \succeq 0} &\ \rho(\bar{P}^{1/2}) \\
\mathrm{s.t.} &\ \sup_{\mathbb{P}_{\xi} \in \mathcal{D}}  \mathbb{P}_{\xi} \left \{ \xi^\top W^\top \bar{P} W\xi> 1 \right \} \le \varepsilon
\end{split}
\label{eq:DRprobdesign01}
\tag{DRFD'}
\end{equation}

\begin{thm}
Solving the following problem always yields a feasible solution to the design problem \eqref{eq:DRprobdesign01} with $\mathcal{D} = \mathcal{D}_{\alpha}(\gamma_1,\gamma_2,\Xi_b)$:
\begin{equation}
\begin{split}
\max_{\bar{P},Q,q,q_0,\eta,\beta,\tilde{\beta}} & ~ \rho(\bar{P}^{1/2}) \\
\mathrm{s.t.}~~~~~ &\ \gamma_2 {\rm Tr} \left \{
QS_0^\alpha \right \} + q_0 \le \varepsilon \eta \\
&~ \eta \ge 0,~\beta_j \ge 0,~\tilde{\beta}_j \ge 0,~j=1,\cdots,N_e \\
&~ \begin{split}
& \begin{bmatrix}
Q - W^\top \bar{P} W & q & 0 \\
q^\top & q_0 - \eta & 0 \\
0 & 0 & 1 \\
\end{bmatrix} + \sum_{j=1}^{N_e} \beta_j \Phi_j \\
& ~~~~~~~~~~~~~~~~~+ {\eta \over \tau_0^\alpha} \begin{bmatrix}
0 & 0 & 0 \\
0 & \alpha + 1 & {-\alpha \over 2\tau_0} \\
0 & {-\alpha \over 2\tau_0} & 0 \\
\end{bmatrix} \succeq 0
\end{split} \\
&~ \begin{bmatrix}
Q & q \\
q^\top & q_0
\end{bmatrix} + \sum_{j=1}^{N_e} \tilde{\beta}_j \tilde{\Phi}_j \succeq 0,~ Q \succeq 0,~ \bar{P} \succeq 0
\end{split}
\label{eq:probdesign17}
\end{equation}
where $\tau_0 > 0$.
\label{thm:10}
\end{thm}
\begin{proof}
In virtue of Theorem \ref{thm:8}, we arrive at the following \textit{inner approximation} to the design problem \eqref{eq:DRprobdesign01}:
\begin{equation*}
\begin{split}
\min_{\bar{P},Q,q,q_0,\eta,\beta,\tilde{\beta},M} & ~ \rho(\bar{P}^{1/2}) \\
\mathrm{s.t.}~~~~~~~ &\ \gamma_2 {\rm Tr} \left \{
QS_0^\alpha \right \} + q_0 \le \varepsilon \\
&~ {\rm Constraints}~\eqref{eq:probdesign16b}-\eqref{eq:probdesign16e} \\
&~ M = W^\top \bar{P} W,~ \bar{P} \succeq 0
\end{split}
\end{equation*}
Because $\{ \bar{P}, M\}$ have to be optimized in conjunction with $\eta$, the constraint \eqref{eq:probdesign16c} now becomes a bilinear matrix inequality (BMI). By defining new variables $Q := Q / \eta$, $q := q / \eta$, $q_0 := q_0 / \eta$, $\beta_j := \beta_j / \eta$, $\tilde{\beta}_j := \tilde{\beta}_j / \eta$, $\eta := 1 / \eta$, bilinear terms can be eliminated, giving rise to the reformulation \eqref{eq:probdesign17}. 
\end{proof}
\begin{rem}
Because both metrics $\rho_1(\bar{P}^{1/2})$ and $\rho_2(\bar{P}^{1/2})$ are concave in $\bar{P}$ \cite{boyd2004convex}, the design problem \eqref{eq:probdesign17} is always a convex program amenable to off-the-shelf solvers.
\end{rem}
\begin{rem}
As an approximation, the pessimisim underlying Theorem \ref{thm:10} is closely related to the selection of $\tau_0 > 0$. In order for a good approximation, gridding of $\tau_0$ can be carried out around \eqref{eq:optTau0} efficiently and the best solution can be chosen from all suboptimal candidates. Given the optimal solution $\bar{P}^*$ to \eqref{eq:probdesign17}, the projection matrix $P^*$ can be obtained from the Cholesky decomposition $\bar{P}^* = LL^\top$.
\end{rem}

\subsection{Safe thresholding for change detection}
Indeed, the usage of developed generalized Gauss bounds is not limited to integrated fault detection design. We can envisage its broader applicability in versatile change detection tasks across different thematic fields of systems and control, where \textit{anomaly detectors in quadratic forms} ${\rm Index}(\xi) = \xi^\top M \xi$ play a central role. Some notable instances are described below.
\begin{itemize}
  \item \textbf{Attack detection in CPS}. Very often, CPS is modelled as a discrete-time linear time-invariant system. Using a state estimator, one defines residuals as output errors, based on which a quadratic failure detector is commonly constructed \cite{mo2015performance,renganathan2020distributionally}.
  \item \textbf{Control performance monitoring (CPM)}. To detect performance degradation, the celebrated Harris index in minimum variance control and quadratic cost in linear quadratic Gaussian (LQG) control are widely adopted as CPM indices, see e.g. \cite{huang2008dynamic,ding2021control}.
  \item \textbf{Multivariate statistical process monitoring}. Popular monitoring indices such as the Hotelling's $T^2$ and the squared prediction error can be constructed with multivariate analysis methods including principal component analysis and slow feature analysis \cite{joe2003statistical,shang2015concurrent}.  
\end{itemize}

Given an index ${\rm Index}(\xi)$, its realistic FAR is dependent on the alarm threshold $J_{\rm th}$, which is generically calibrated under the Gaussian assumption on $\mathbb{P}_{\xi}$. A common concern in above thematic fields is that the Gaussianity itself could be unjustified in engineering practice, and thus robustness against ambiguity and/or variations in probability distributions is desired. In this case, $\mathcal{D}_\alpha$ acts as an appealing choice to delineate distribution of $\xi$, and in order to reliably calibrate $J_{\rm th}$ under a given confidence, generalized Gauss bounds in Theorems \ref{thm:4} and \ref{thm:8} provide a practically useful alternative to generic $\chi^2$-quantiles in removing false alarms. Here we only discuss the case of $\alpha$-unimodal set $\mathcal{D}_\alpha(\gamma_1,\gamma_2,\Xi_b)$ with bounded support, by making use of Theorem \ref{thm:8}. The proof is omitted since it requires no new ideas.

\begin{thm} Given a quadratic anomaly detector ${\rm Index}(\xi) = \xi^\top M \xi$, where the distribution of $\xi$ is included in the ambiguity set $\mathcal{D} = \mathcal{D}_\alpha(\gamma_1,\gamma_2,\Xi_b)$ and $\Xi_b$ is given in \eqref{eq:support}. An alarm threshold $J_{\rm th}$ rendering FAR no higher than the tolerance $\varepsilon \in (0,1)$ can be obtained by solving the following SDP:
\begin{equation*}
\begin{split}
\min_{Q,q,q_0,\eta,\beta,\tilde{\beta},J_{\rm th}} & ~J_{\rm th} \\
\mathrm{s.t.}~~~~~~ & ~ \gamma_2 {\rm Tr} \left \{
QS_0^\alpha \right \} + q_0 \le \varepsilon J_{\rm th} \\
&~ \eta \ge 0,~\beta_j \ge 0,~\tilde{\beta}_j \ge 0,~j=1,\cdots,N_e \\
& \begin{bmatrix}
Q - \eta M & q & 0 \\
q^\top & q_0 - J_{\rm th} & 0 \\
0 & 0 & \eta \\
\end{bmatrix} + \sum_{j=1}^{N_e} \beta_j \Phi_j  \\
& ~~~~~~~~~~~~~~~~ +{J_{\rm th} \over \tau_0^{\alpha}} \begin{bmatrix}
0 & 0 & 0 \\
0 & \alpha + 1 & {-\alpha \over 2\tau_0} \\
0 & {-\alpha \over 2\tau_0} & 0 \\
\end{bmatrix}  \succeq 0 \\
&~ \begin{bmatrix}
Q & q \\
q^\top & q_0
\end{bmatrix} + \sum_{j=1}^{N_e} \tilde{\beta}_j \tilde{\Phi}_j \succeq 0, ~ Q \succeq 0
\end{split}
\end{equation*}
where $\tau_0 > 0$ is a tuning parameter.
\label{thm:11}
\end{thm}

\section{Case Studies}
\begin{figure}[ht]
  \centering
  \includegraphics[width=0.4\textwidth]{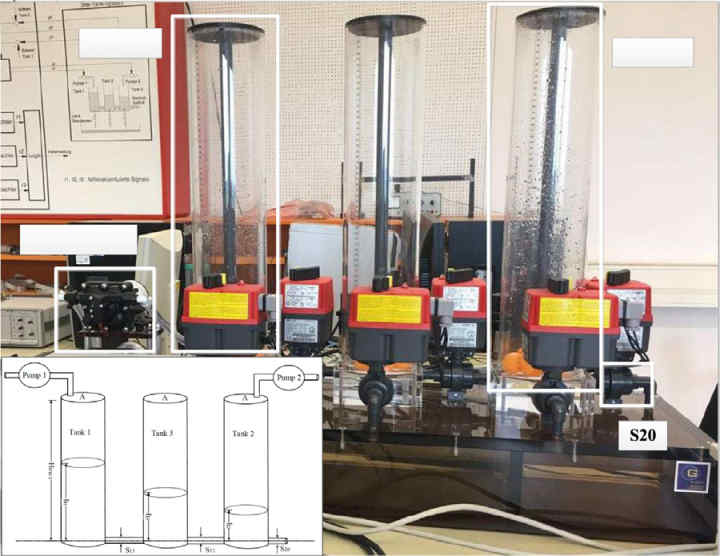}\\
  \caption{Schematic of the experimental three-tank system TTS20.}\label{fig:threetank}
\end{figure}

Next, we perform case studies using realistic data collected from an experimental three-tank system. A pictorial description of this apparatus is given in Fig. \ref{fig:threetank}, where three tanks are interconnected via two pipelines. Water is fed into Tanks 1 and 2 by two pumps, whose flowrates act as system inputs $u(k) \in \mathbb{R}^2$. System states $x(k) = [x_1(k) ~ x_2(k) ~ x_3(k)]^\top$ are the levels of three tanks, among which $y(k) = [x_1(k) ~ x_3(k)]^\top$ constitute the outputs. All system matrices in state-space equations are derived by linearizing differential equations around the operating point with a sampling interval $\Delta t = 5$s \cite{shang2021distributionally}. Then, the parity space approach \cite{chow1984analytical} is adopted to construct a residual signal $v(k) \in \mathbb{R}^9$ with $s = 6$, whose dynamics is governed by $v(k) = W_v d_s(k) + V_v f_s(k)$.

Following Remark \ref{rem:1}, we seek to obtain a ``calibrated" residual $r(k) = Pv(k)$ with distributional robustness, by optimizing $P$ while regarding the fault-free realization $v_0(k) = W_v d_s(k)$ as uncertainty $\xi$ that follows an unknown distribution. Data samples of $v_0(k)$ are collected in routine fault-free operations, based on which statistical information of $\xi$ such as covariance and support can be estimated. In particular, the support $\Xi_b$ is determined as a hyper-rectangular $\Xi_b = \{ \xi | | \xi_i | \le 1.2 \times \xi_i^{\max},~\forall i =1,\cdots,9\}$ that reliably covers all samples and is also representable as \eqref{eq:support}, where $\xi_i^{\max}$ is the element-wise maximum of $|\xi_i|$ on empirical samples. Meanwhile, unimodality can be inspected from scatter plots of $v_0(k)$; hence, we simply choose $\alpha = n = 9$, which corresponds to star-unimodality, to encode such minimal structral information. Using these information, various ambiguity sets can be constructed, and for a particular metric $\rho_i(\cdot),~i=1,2$, the following design schemes are developed.
\begin{itemize}
\item DR$i$-U: Optimal solution to \eqref{eq:DRprobdesign0} under $\mathcal{D}(\gamma_1,\gamma_2,\mathbb{R}^n)$ assuming unbounded support \cite[Theorems 2 and 3]{shang2021distributionally};
\item DR$i$-U-$\alpha$: Feasible solution to \eqref{eq:DRprobdesign0} under $\mathcal{D}_\alpha(\gamma_1,\gamma_2,\mathbb{R}^n)$ assuming unbounded support and $\alpha$-unimodality in Theorems \ref{thm:6} and \ref{thm:7};
\item DR$i$-B: Feasible solution to \eqref{eq:DRprobdesign0} under $\mathcal{D}(\gamma_1,\gamma_2,\Xi_b)$ assuming bounded support \cite[Theorem 4]{shang2021distributionally};
\item DR$i$-B-$\alpha$: Feasible solution to \eqref{eq:DRprobdesign0} under $\mathcal{D}_\alpha(\gamma_1,\gamma_2,\Xi_b)$ assuming bounded support and $\alpha$-unimodality in Theorem \ref{thm:10}, with $\tau_0$ optimally calibrated by gridding;
\end{itemize}

\begin{figure}[ht]
  \centering
  \includegraphics[width=0.48\textwidth]{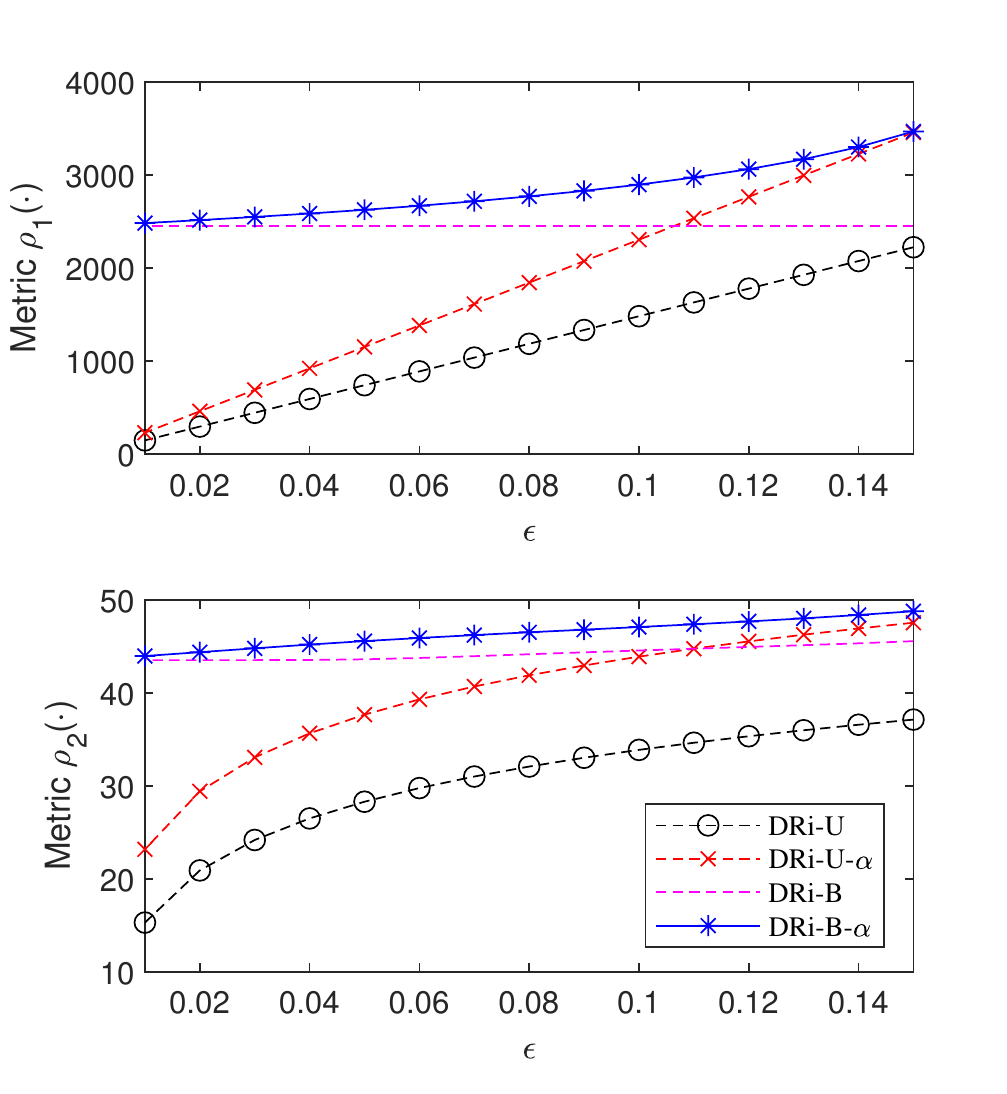}\\
  \caption{Optimal values of problem \eqref{eq:DRprobdesign0} under different ambiguity sets and detectability metrics $\rho(\cdot)$.}\label{fig:obj}
\end{figure}

We first investigate the effect of introducing $\alpha$-unimodality in reducing conservatism. The problem \eqref{eq:DRprobdesign0} is resolved under various choices of ambiguity sets and different tolerances $\varepsilon$, and all induced convex programs are resolved with \texttt{mosek}. The resulting optimal values of \eqref{eq:DRprobdesign0} are displayed in Fig. \ref{fig:obj}. It can be seen that by introducing unimodality information, higher values of detectability metrics are attained, indicating a smaller feasible region and thus a reduction of pessimism. Under a relatively large $\varepsilon$, the effect of assuming unimodality tends to outweigh the usage of bounded support information. When $\varepsilon$ is small, DR$i$-B outdoes DR$i$-U because DR$i$-B tends to classic set-membership robust design, while the worst-case distribution in DR$i$-U has masses outside the support and is thus unrealistic. By synthesizing all information, the best detectability can always be attained by DR$i$-B-$\alpha$.

\begin{figure*}[t]
  \centering
  \includegraphics[width=0.7\textwidth]{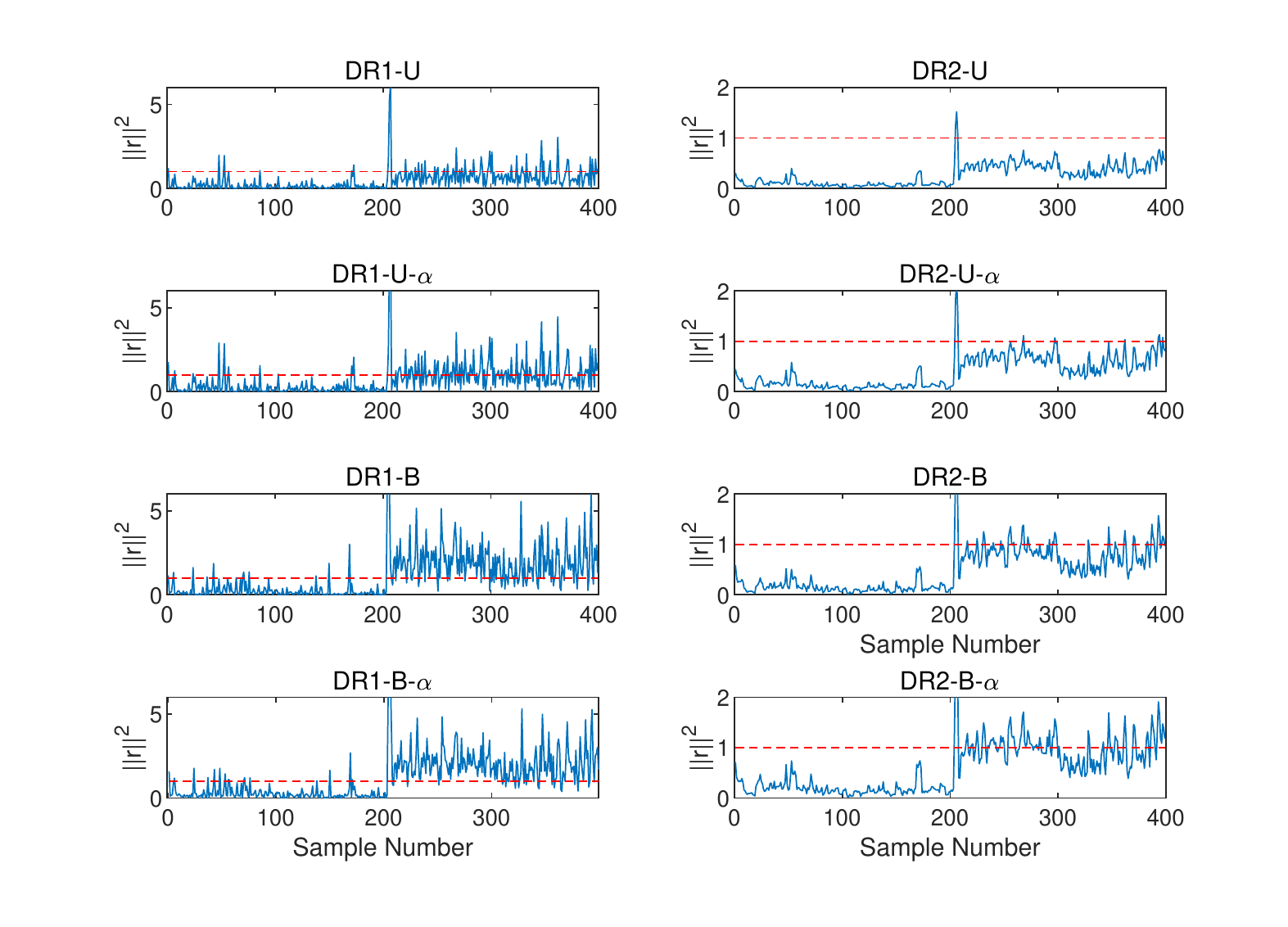}\\
  \caption{Detection of leakage in Tank 3 using different design matrix. Red dashed lines denote detection thresholds.}\label{fig:residual}
\end{figure*}

Then we investigate detection performance using a test dataset collected from the apparatus, which consists of both normal and faulty samples. The first 200 samples are collected under healthy conditions, based on which FAR can be evaluated, while a fault of leakage in Tank 3 is introduced from the 201st sample till the end, based on which FDR can be evaluated. Using combinations of different metric $\rho_i(\cdot)$ and ambiguity set, their detection results are depicted in Fig. \ref{fig:residual}, and performance indices are summarized in Tables \ref{tab:stat_rho1} and \ref{tab:stat_rho2}. For $\rho_1(\cdot)$, all induced detectors successfully keep FARs lower than the tolerance, showcasing the guaranteed robustness. In contrast, the classic GLRT regime assuming Gaussianity yields an FAR of $45\%$, resulting in massive nuisance alarms. On the other hand, the detectability of both DR1-U and DR1-B is improved by their $\alpha$-unimodal counterparts, where a desirable tradeoff is achieved by DR1-B-$\alpha$ with the highest FDR $87\%$. In general, the DRFD designs induced by $\rho_2(\cdot)$ appear to be more conservative. In DR2-U, DR2-U-$\alpha$, and DR2-B, the fault remains largely undetected, while the exploitation of both support and unimodality information in DR2-B-$\alpha$ helps better showcase the fault and reduce the pessimism.

\begin{table}[htp]
\caption{Fault Detection Performance with Metric $\rho_1(\cdot)$}
\centering
\begin{tabular}{c c c c c}
\toprule
    & DR1-U & DR1-U-$\alpha$ & DR1-B & DR1-B-$\alpha$ \\
\midrule
FAR (\%) & 3.50 & 6.00 & 6.00 & 7.50 \\
FDR (\%) & 24.50 & 46.00 & 76.50 & 87.00 \\
\bottomrule
\end{tabular}
\label{tab:stat_rho1}
\end{table}

\begin{table}[htp]
\caption{Fault Detection Performance with Metric $\rho_2(\cdot)$}
\centering
\begin{tabular}{c c c c c}
\toprule
    & DR2-U & DR2-U-$\alpha$ & DR2-B & DR2-B-$\alpha$ \\
\midrule
FAR (\%) & 0.00 & 0.00 & 0.00 & 0.00 \\
FDR (\%) & 1.50 & 4.50 & 21.50 & 46.00 \\
\bottomrule
\end{tabular}
\label{tab:stat_rho2}
\end{table}

\section{Concluding Remarks}
In this article, a new distributionally robust design scheme is developed to maximize fault detectability and regulate false alarms without requiring precise distribution knowledge. By further integrating $\alpha$-unimodality information, the conservatism of previous moment-based DRFD schemes can be reliably alleviated. To tackle the worst-case constraint on FAR, we first establish a new multivariate $\alpha$-unimodal Gauss bound on the tail probability outside an ellipsoid, which strictly improves upon its Chebyshev counterpart with the improvement factor explicitly given. Based on this, feasible solutions to DRFD problems are derived in closed form. Then we develop a tightened $\alpha$-unimodal Gauss bound by further injecting support information, which enables us to tackle the DRFD design problem approximately by solving a convex program and further alleviate the design conservatism. The developed multivariate Gauss bounds also apply to a broader class of change detection tasks across different areas in systems and control. The efficacy of the new fault detection design scheme in reducing conservatism is illustrated using data collected from an experimental three-tank apparatus.

\section*{Appendix: Proof of Theorem \ref{thm:4}}
We first derive a convex approximation of \eqref{eq:5}. Define $h(\tau) = - \tau^{-\alpha}$, which is a concave in $\tau > 0$. Then given $\tau_0 := \| M^\frac12 \xi_0 \|$ and $\tau := \| M^\frac12 \xi \| > 0$, $h(\tau)$ is majorized by its first-order Taylor approximation:
\begin{equation*}
h(\tau) \le h(\tau_0) + (\tau - \tau_0)  h'(\tau_0) = -\tau_0^{-\alpha} + \alpha(\tau - \tau_0)
\end{equation*}
which amounts to
\begin{equation*}
\begin{split}
-\| M^\frac12 \xi \|^{-\alpha} \le -(\alpha + 1)\tau_0^{-\alpha} + \alpha \tau_0^{-\alpha - 1} \| M^\frac12 \xi \|.
\end{split}
\end{equation*}
One then obtains an inner approximation to \eqref{eq:5}:
\begin{subequations} \label{eq:thm1_apprx}
\begin{align}
\min_{Q,q,q_0} &\ \gamma_2 {\rm Tr} \left \{
QS_0^\alpha \right \} + q_0  \\
\mathrm{s.t.}~ &~ \begin{bmatrix}
Q & q \\
q^\top & q_0
\end{bmatrix} \succeq 0 \label{eq:thm1_apprx_b} \\
\begin{split} & ~~\xi^\top Q\xi + 2\xi^\top q + q_0 - 1 + (\alpha + 1)\tau_0^{-\alpha} \\
&~~~~~~~~~~~~~~~~~\ge \alpha \| M^\frac12 \xi \| \cdot \tau_0^{-\alpha - 1},~ \forall \xi \in \mathbb{R}^n \end{split} \label{eq:thm1_apprx_c}
\end{align}
\end{subequations}
Note that $Q \succeq 0$ is implied by \eqref{eq:thm1_apprx_b} and is thus omitted. Next we focus on the semi-infinite constraint \eqref{eq:thm1_apprx_c}, which can be restated as that $\forall \xi, t$ such that $\xi^\top M\xi \ge t^2$, the inequality
\begin{equation*}
 \xi^\top Q\xi + 2\xi^\top q + q_0 - 1 + (\alpha + 1)\tau_0^{-\alpha}\ge \alpha \tau_0^{-\alpha - 1} t
\end{equation*}
always holds. In virtue of the S-procedure \cite{yakubovich1997s}, the semi-infinite constraint equals to the existence of Lagrangian multiplier $\eta \ge 0$ such that
\begin{equation}
g(\xi,t) \ge 0,~\forall \xi \in \mathbb{R}^n,~ t \in \mathbb{R}
\label{eq:thm1_semiinf}
\end{equation}
where
\begin{equation*}
\begin{split}
&~ g(\xi,t)  \\
\triangleq &~ \xi^\top Q\xi + 2\xi^\top q + q_0 - 1 + (\alpha + 1)\tau_0^{-\alpha} - \alpha \tau_0^{-\alpha - 1} t + \eta t^2  \\
&~  - \eta\xi^\top M \xi
\end{split}
\end{equation*}
Note that \eqref{eq:thm1_semiinf} can be rewritten as an LMI. Consequently, problem \eqref{eq:thm1_apprx} equals to the following SDP
\begin{subequations} \label{eq:thm1_apprx1}
\begin{align}
\min_{Q,q,q_0,\eta} &\  \gamma_2 {\rm Tr} \left \{
QS_0^\alpha \right \} + q_0  \\
\mathrm{s.t.}~~ & \begin{bmatrix}
Q & q \\
q^\top & q_0
\end{bmatrix} \succeq 0,~ \eta \ge 0 \label{eq:thm1_apprx1b} \\
\begin{split} & \begin{bmatrix}
Q - \eta M & q & 0 \\
q^\top & q_0-1 & 0 \\
0 & 0 & \eta \\
\end{bmatrix} + {1 \over \tau_0^{\alpha}} \begin{bmatrix}
0 & 0 & 0 \\
0 & \alpha + 1 & {-\alpha \over 2\tau_0} \\
0 & {-\alpha \over 2\tau_0} & 0 \\
\end{bmatrix} \succeq 0
\end{split} \label{eq:thm1_apprx1c}
\end{align}
\end{subequations}
where \eqref{eq:thm1_semiinf} is recast as \eqref{eq:thm1_apprx1c}. Since $\tau_0 \ge 1$, it holds that $\alpha\tau_0^{-\alpha - 1}/2 > 0$, and with the aim of securing the positive semi-definiteness of \eqref{eq:thm1_apprx1c}, it must be that $\eta > 0$ and $q_0 + (\alpha +1) \tau_0^{-\alpha} - 1 > 0$. As a result, according to the Schur's complement, the LMI \eqref{eq:thm1_apprx1c} can be recast as:
\begin{subequations} \label{eq:thm1_apprx2}
\begin{align}
& Q - \eta M \succeq 0 \label{eq:thm1_apprx2a} \\
& Q - \eta M \succeq \frac{qq^\top}{q_0 + (\alpha +1) \tau_0^{-\alpha} - 1} \label{eq:thm1_apprx2b} \\
& \eta \left [ q_0 + (\alpha +1) \tau_0^{-\alpha} - 1 \right ] \ge \frac{\alpha^2}{4}\tau_0^{-2\alpha - 2} \label{eq:thm1_apprx2c} \\
& \eta > 0,~ q_0 + (\alpha +1) \tau_0^{-\alpha} - 1 > 0 \label{eq:thm1_apprx2d}
\end{align}
\end{subequations}
where \eqref{eq:thm1_apprx2a} is implied by \eqref{eq:thm1_apprx2b} and thus can be omitted. Next we distinguish between three cases.

\textit{Case 1:} $(\alpha +1) \tau_0^{-\alpha} \le 1$. In this case, $q_0 > 1 - (\alpha +1) \tau_0^{-\alpha} \ge 0$. It follows from the Schur's complement that the LMI \eqref{eq:thm1_apprx1b} reduces to $Q \succeq qq^\top / q_0$. Thus problem \eqref{eq:thm1_apprx1} becomes:
\begin{equation*}
\begin{split}
\min_{Q,q,q_0,\eta} &\ \gamma_2 {\rm Tr} \left \{
QS_0^\alpha \right \} + q_0  \\
\mathrm{s.t.}~~ &~ Q \succeq {qq^\top \over q_0},~Q - \eta M \succeq \frac{qq^\top}{q_0 + (\alpha +1) \tau_0^{-\alpha} - 1} \\
                       &~ \eta \left [ q_0 + (\alpha +1) \tau_0^{-\alpha} - 1 \right ] \ge \frac{\alpha^2}{4}\tau_0^{-2\alpha - 2} \\
                       &~ q_0 + (\alpha +1) \tau_0^{-\alpha} - 1 > 0
\end{split}
\end{equation*}
which amounts to
\begin{equation*}
\begin{split}
\min_{q,q_0,\eta,r} &~ r + q_0  \\
\mathrm{s.t.}~ & ~ r \ge {\gamma_2 q^\top S_0^\alpha q \over q_0} \\
                        & ~ r \ge \frac{\gamma_2 q^\top S_0^\alpha q}{q_0 + (\alpha +1) \tau_0^{-\alpha} - 1} + \gamma_2 \eta {\rm Tr}\{ MS_0^\alpha \} \\
                       &~ \eta \left [ q_0 + (\alpha +1) \tau_0^{-\alpha} - 1 \right ] \ge \frac{\alpha^2}{4}\tau_0^{-2\alpha - 2} \\
                       &~ q_0 + (\alpha +1) \tau_0^{-\alpha} - 1 > 0
\end{split}
\end{equation*}
Clearly, it holds that $q^* = 0$, $r^* = \gamma_2 \eta {\rm Tr}\{ MS_0^\alpha \}$, which yields:
\begin{equation}
\begin{split}
& \left \{ \begin{split}
\min_{q_0,\eta} &\  \gamma_2 \eta {\rm Tr}\{ MS_0^\alpha  \} + q_0  \\
\mathrm{s.t.} &~ \eta \left [ q_0 + (\alpha +1) \tau_0^{-\alpha} - 1 \right ] \ge \frac{\alpha^2}{4}\tau_0^{-2\alpha - 2} \\
                       &~ q_0 + (\alpha +1) \tau_0^{-\alpha} - 1 > 0
\end{split} \right . \\
= & \left \{ \begin{split}
\min_{q_0} &\ \frac{\alpha^2\tau_0^{-2\alpha - 2} \gamma_2 {\rm Tr}\{ MS_0^\alpha \}}{4[q_0 + (\alpha +1) \tau_0^{-\alpha} - 1]} + q_0 \\
\mathrm{s.t.}  &~ q_0 + (\alpha +1) \tau_0^{-\alpha} - 1 > 0
\end{split} \right . \\
= &~ \alpha \tau_0^{-\alpha - 1}\sqrt{\gamma_2{\rm Tr}\{ MS_0^\alpha \}} - (\alpha + 1)\tau_0^{-\alpha} + 1 \\
\triangleq &~ f_1(\tau_0)
\end{split}
\label{eq:thm1_result1}
\end{equation}
where
\begin{equation*}
\begin{split}
q_0^* & = {\alpha \over 2} \tau_0^{-\alpha - 1}\sqrt{\gamma_2{\rm Tr}\{ MS_0^\alpha \}} + 1 - (\alpha +1) \tau_0^{-\alpha} \\
& > 1 - (\alpha +1) \tau_0^{-\alpha}
\end{split}
\end{equation*}
is attainable.

\textit{Case 2:} $1 < (\alpha +1) \tau_0^{-\alpha} \le 1 + {\alpha \over 2} \tau_0^{-\alpha - 1}\sqrt{\gamma_2{\rm Tr}\{ QS_0^\alpha \}}$. If $q_0 = 0$, \eqref{eq:thm1_apprx1b} boils down to $Q \succeq 0$ with $q^* = 0$, and consequently the optimal value of \eqref{eq:thm1_apprx1} is obtained as
\begin{equation}
f_2(\tau_0) \triangleq \frac{\alpha^2\tau_0^{-2\alpha - 2} \gamma_2 {\rm Tr}\{ MS_0^\alpha \}}{4[(\alpha +1) \tau_0^{-\alpha} - 1]}.
\label{eq:17}
\end{equation}
If $q_0 > 0$ otherwise, proceeding in a similar manner to Case 1, one obtains the optimal value as $f_1(\tau_0)$ in \eqref{eq:thm1_result1}. Because the value of \eqref{eq:17} is lower than that of \eqref{eq:thm1_result1}, one concludes that the optimal value is $f_1(\tau_0)$, and thus Case 2 can be combined with Case 1 as a single case.

\textit{Case 3:} $(\alpha +1) \tau_0^{-\alpha} > 1 + {\alpha \over 2} \tau_0^{-\alpha - 1} \sqrt{\gamma_2{\rm Tr}\{ QS_0^\alpha \}}$. It is an easy exercise to verify that $q_0^* = 0$, and thus the optimal value is $f_2(\tau_0)$.

Summarizing above cases yields the optimal value of the relaxed problem \eqref{eq:thm1_apprx1}:
\begin{equation}
\begin{split}
&~ \sup_{\mathbb{P} \in \mathcal{D}} \mathbb{P}_{\xi}  \left \{ \xi^\top M \xi > 1 \right \}  \\
\le & \left \{
\begin{split}
& f_1(\tau_0),~{\rm if}~(\alpha +1) \tau_0^{-\alpha} \le 1 + {\alpha \over 2} \tau_0^{-\alpha - 1}\sqrt{\gamma_2{\rm Tr}\{ MS_0^\alpha \}} \\
& f_2(\tau_0),~{\rm otherwise}
\end{split}
\right .
\end{split}
\label{eq:18}
\end{equation}

Next, we seek the best approximation among all choices of $\tau_0 > 0$. For convenience we define $z = 1/\tau_0 \in (0,1]$, $g_1(z) = f_1(\tau_0)$, $g_2(z) = f_2(\tau_0)$ and $g(z) = 1 + {\alpha \over 2} \sqrt{\gamma_2{\rm Tr}\{ MS_0^\alpha \}} z^{\alpha + 1} - (\alpha + 1)z^\alpha$. In this way, \eqref{eq:18} becomes:
\begin{equation*}
\begin{split}
&~ \sup_{\mathbb{P} \in \mathcal{D}} \mathbb{P}_{\xi}  \left \{ \xi^\top M\xi > 1 \right \} \le h(z) \triangleq \left \{
\begin{split}
& g_1(z),~{\rm if}~g(z) \ge 0 \\
& g_2(z),~{\rm if}~g(z) < 0
\end{split}
\right .
\end{split}
\end{equation*}
where $h(z)$ is continuous. Thus, it suffices to resolve $\min_{z \in (0,1]} h(z)$. By setting $g_i'(z) = 0~(i=1,2)$, one obtains a unique stationary point of $g_i(z)$ on $\mathbb{R}_+$:
\begin{equation*}
z_1^* = \frac{1}{\sqrt{\gamma_2{\rm Tr} \{ MS_0^\alpha \}}},~ z_2^* = \sqrt{c_\alpha}.
\end{equation*}
Note that $0<z_2^*<1$ always holds. Meanwhile, $g_i(z)$ is first decreasing on $(0,z_i^*)$ and then increasing. As for $g(z)$, there is also a unique stationary point $z^* = 2 / \sqrt{\gamma_2{\rm Tr} \{ MS_0^\alpha \}} > 0$. Next, the following cases are distinguished.

\textit{Case 1:} $g(z^*) \ge 0$, which equals to $\gamma_2 {\rm Tr} \{ MS_0^\alpha \} \ge 4$. In this case, $g(z) \ge 0$ always holds for $z \in (0,1]$, and the minimizer $z_1^* \le 1/2$ of $g_1(z)$ is always attainable. Thus,
\begin{equation*}
\min_{z \in (0,1]} h(z) = g_1(z_1^*) = 1 - \frac{1}{(\gamma_2{\rm Tr} \{ MS_0^\alpha \})^{\alpha / 2}}.
\end{equation*}

\textit{Case 2:} $g(z^*) < 0$ and $g(z_2^*) \ge 0$, which is equivalent to:
\begin{equation*}
{1 \over c_\alpha} \le \gamma_2{\rm Tr} \{ MS_0^\alpha \} < 4,
\end{equation*}
thereby indicating $0 < z_1^* \le z_2^* < 1$. Moreover, it can be verified that $g(1) < 0$. Thus, $g(z)$ has a unique root $\hat{z} \in (z_2^*, 1)$, and one obtains:
\begin{equation}
h(z) = \left \{
\begin{split}
& g_1(z),~{\rm if}~z \in (0,\hat{z}] \\
& g_2(z),~{\rm if}~z \in (\hat{z},1]
\end{split}
\right .
\label{eq:23}
\end{equation}
Note that $\min_{z \in (0,\hat{z}]} g_1(z) = g_1(z_1^*)$, while $\min_{z \in (\hat{z},1]} g_2(z) = g_1(\hat{z}) \ge g_1(z_1^*)$. It immediately follows that $\min_{z\in(0,1]} h(z) = g_1(z_1^*)$.

\textit{Case 3:} $g(z^*) < 0$ and $g(z_2^*) < 0$. In this case, it holds that
\begin{equation*}
\gamma_2 {\rm Tr} \{ MS_0^\alpha \} < {1 \over c_\alpha},
\end{equation*}
which implies $g(z_1^*) < 0$ and $z_1^* > z_2^*$. Because there always exists $\hat{z} \in (0, z_2^*)$ such that $g(\hat{z}) = 0$, $h(z)$ can be expressed as \eqref{eq:23}, and we have $\min_{z \in (0,\hat{z}]} g_1(z) = g_2(\hat{z}) \ge g_2(z_2^*)$ and $\min_{z \in (\hat{z},1]} g_2(z) = g_2(z_2^*)$. This gives rise to
\begin{equation*}
\min_{z \in (0,1]} h(z) = g_2(z_2^*) = c_\alpha \gamma_2 {\rm Tr} \{ MS_0^\alpha \}.
\end{equation*}

Summarizing above cases yields \eqref{eq:10}.


\end{document}